\documentclass[12pt]{article}

\usepackage[a4paper,margin=1in]{geometry}
\usepackage{setspace}
\usepackage[english]{babel}
\usepackage{microtype}

\usepackage{amsmath, amssymb, amsthm, mathtools, mathrsfs}

\usepackage{eufrak}

\usepackage{graphicx}
\usepackage{float}
\usepackage{caption}
\usepackage{subcaption}

\usepackage{hyperref}
\usepackage{cleveref}
\usepackage{url}

\usepackage{enumitem}

\newcounter{cases}
\newcounter{subcases}[cases]

\title{\bfseries Delayed Pattern Formation in Two-Dimensional Domains}

\author{
Nirmali Prabha Das,\thanks{Bolyai Institute, University of Szeged, Hungary}
Istv\'an Bal\'azs,\footnotemark[1]
Bornali Das,\footnotemark[1]
Gergely R\"ost\footnotemark[1]
}

\date{}

\begin{document}
\maketitle

\begin{abstract}
This study investigates how the interaction between gene expression time delay and domain size governs spatio-temporal pattern formation in a reaction-diffusion system. To investigate these phenomena, we utilize a modified version of the Schnakenberg model called the ligand internalisation (LI) model. In a one-dimensional domain, a linear relationship has been observed between the gene expression time delay and the time it takes for patterns to form. We extend the model to the two-dimensional domain and confirm that a similar relationship holds there as well. However, our exploration reveals a non-monotonic correlation between domain size and the time required for pattern emergence.

To unravel these dynamics, we consider a range of initial conditions, including random perturbations of the spatially homogeneous steady state and initial conditions from its unstable manifold. We compute a two-parameter chart of patterns with respect to time delay and domain size.
\end{abstract}


%% \linenumbers

%% main text

\parskip=5pt plus 1pt

\section{Introduction}\label{sec1}
Biological pattern formation is a ubiquitous phenomenon in nature, yet its origin remains a central problem in developmental biology. How a spatially homogeneous population of cells spontaneously breaks symmetry to generate ordered structures is a question that has intrigued researchers for decades. The spontaneous emergence of spatial organization from an initially disordered system has motivated extensive theoretical and experimental investigations across biology, chemistry, physics, and ecology.

The development of a spatial pattern of tissue structures is one of the most fundamental morphogenesis processes. Following a set of intricate phenomena, an organism evolves from an almost homogeneous tissue. It has been established using a mathematical model of remarkable simplicity founded on auto- and cross-catalysis that such a basic molecular process can account for the foundational characteristics of pattern formation \cite{gierer1972theory}. 
The application of reaction-diffusion systems to study spatial pattern  formation has a long history within the broader framework of self-organizing systems, beginning with Alan Turing's publication in 1952 
\cite{turing1952chemical}. Reaction-diffusion systems have since become a paradigm of pattern development in chemistry and biology, wherein such an order was established as the governing principle of biological morphogenesis, the mechanism by which organisms undergo metamorphism. The terminology "morphogen" was first introduced in this context by Alan Turing to designate chemical substances whose spatial distribution is the physical foundation of the underlying principle, as cells react to variations in concentration. Alan Turing postulated that the interaction of morphogens, the evolution of which is defined by coupled reaction-diffusion equations, could be used to quantitatively model the pattern generation process on a strictly chemical basis. Perturbation from the steady state of a particular chemical system results in the emergence of structure, from an initially homogeneous state. The notion of morphogens has since served as a classical framework to examine and comprehend various aspects of plant and animal development including microbiological synthesis processes, pigmentation patterns, ecological dispersion, and chaos, as well as disease transmission dynamics \cite{rybakin2001morphogenesis, green2002morphogen, grieneisen2012morphogengineering, murray2001mathematical}.

The fundamental processes of all living organisms such as reproduction, growth and development, adaptation, homeostasis, evolution, or environmental responsiveness are regulated by gene expression.  The reaction-diffusion equations, which form a system of partial differential equations, have been extensively researched and are widely utilized to generate spatially diverse patterned states from homogeneous symmetry breaking via the Turing instability \cite{meinhardt1987model, bard1981model}. Although there are many prototype "Turing systems" accessible, it can be challenging to determine their characteristics, functional forms, and overall suitability for a particular application.

The study of Turing pattern formation for particular reaction-diffusion systems that have been proposed to model different morphogenetic processes has accumulated a substantial volume of scientific literature \cite{grieneisen2012morphogengineering}. Various applications can be found that have been used to study the coat pattern formation \cite{sander2003pattern} and pigmentation in animal and plant kingdoms:  the origin of zebra stripes \cite{bard1981model, murray1981pattern, murray1981pre}, sea shell markings \cite{meinhardt1987model}, the scaled pattern on reptiles \cite{cocho1987discrete}, wing patterns of butterflies \cite{nijhout1978wing}, spatial pigment pattern on marine fish \cite{aragon1998spatial}, and pattern formation in plants \cite{meinhardt1998models}. Furthermore, reaction-diffusion systems have been extensively employed in theoretical studies to investigate the effects of a wide range of factors, including but not limited to: spatial inhomogeneity \cite{benson1993diffusion}, growth and curvature effects \cite{plaza2004effect}, boundary conditions \cite{aragon1998spatial}, external stimuli \cite{prigogine1968symmetry},  nonlinear interactions \cite{meinhardt1987model}, temporal dynamics \cite{izhikevich2006fitzhugh},   and evolving dynamic surfaces \cite{kim2020pattern}. It is well-founded that the influence of such factors on the science behind pattern formation is significant and should be thoroughly investigated. Notably, researchers have placed particular emphasis on the role played by domain size and nonlinearities on the instability of solutions and pattern formation. A distinguished example is \cite{kondo1995reaction}, by contributing to the understanding of the underlying mechanism for molecular and cellular interactions that give rise to patterns in developmental biology. By considering Turing patterns in a one-dimensional growing domain they showed the evolution of spatiotemporal patterns that are observed in a certain family of vertebrates (\textit{pomacanthus imperator}, commonly known as the emperor angelfish)  \cite{kondo1995reaction}.

A variety of morphogenesis-related phenomena can be described by considering the interrelationship between two morphogens. Therefore, it is proposed that the embryonic pattern shaped by two morphogens, one playing the role of an activator and the other, an inhibitor, may be a central process in pattern development. Non-linear equations have been presupposed in order to capture the  auto- and cross-catalytic feedback of the activator and inhibitor on their sources \cite{gierer1972theory}.

After the seminal work presented by Turing in 1952,  a significant length of time elapsed before substantial progress was made in developing credible mathematical models to investigate activator-inhibitor dynamics, resulting in the emergence of well-defined stationary and periodic structures. The Gray--Scott model \cite{doelman1997pattern}, the Gierer--Meinhardt model \cite{gierer1972theory}, the FitzHugh--Nagumo Model \cite{toral2003characterization}, the Lengyel--Epstein model \cite{jensen1994localized}, and the Schnakenberg model \cite{schnakenberg1979simple} are some of the well-known models that have been used to study Turing pattern formation.

These models can aid in the understanding of the fundamental mechanism of morphogenesis. Among them, the Schnakenberg model, first developed by Schnakenberg in 1979 is one of the most widely opted chemical reaction-diffusion examples by mathematicians which can be implemented to inspect a wide range of patterns analogous to those found in nature, such as morphogenesis, and can be simply described as: 

\begin{align}
	\begin{aligned}
		\label{eq:L1}
		u_t&=d_u{\nabla^2 u} + f(u,v), \\
		v_t&=d_v{\nabla^2 u} + g(u,v),\\    
	\end{aligned}
\end{align}
where $f(u,v)$ and $g(u,v)$ are the kinetic functions that can be modified to adjust or implement biologically significant attributes. In the case of the Schnakenberg model, $f(u,v) =  a - u + u^2v$, and $g(u,v) =  b - u^2v$, where $u$ and $v$ are respectively the concentrations of the activator and inhibitor. Chemical concentration parameters that control the reaction process are represented by parameters $a$ and $b$ which are positive constants, $d_{u}$ and $d_{v}$ represents the diffusion coefficients of the activator $u$ and inhibitor $v$ respectively. $\nabla^2 u$  and $\nabla^2 v$ are the Laplacian operators that govern the propagation of the morphogens in space.

These earlier studies have typically not considered any time delays in the reaction kinetics, in spite of the fact that the dynamics of pattern generation in reaction-diffusion systems have been found to be significantly impacted by time delays \cite{sargood2022fixed, gaffney2006gene}. Well-founded experimental studies have concluded that the process of intracellular gene expression combines multiple steps including gene transcription and gene translation.  These sub-processes can be time-consuming, with time delays generally being in the range of minutes, but in certain situations, they may last for several hours \cite{tennyson1995human}. The structural resilience of Turing models can significantly vary when temporal delays are included; specifically, pattern formation may have a complex relationship that is highly reliant on the delay. More realistic reaction-diffusion models must account for the interactions at the boundary of the region of self-organizing patterns, its surrounding environment, and any asymmetries there \cite{arcuri1986pattern}.

In this regard, the works of Gaffney and colleagues \cite{sargood2022fixed, gaffney2006gene} that have been designed to study the effect of delay in reaction-diffusion systems are particularly groundbreaking. In one of their earliest works, time delay corresponding to gene expression for protein synthesis has been incorporated into reaction kinetics for one-dimensional pattern-forming Turing systems. The effect of time delay was investigated via linear stability analysis on stationary and steadily growing domains. It has been shown that the delay can significantly lengthen the amount of time needed for the establishment of a stable pattern. Moreover, it can result in a dramatic lack of pattern formation on exponentially growing domains \cite{gaffney2006gene}. 

In the framework of reaction-diffusion systems integrating gene expression time delays to study  Turing pattern formation, the activator-depleted model, commonly referred to as the Schnakenberg model \cite{yi2017bifurcation} has received a great deal of attention. The ligand-internalisation (LI) and reverse ligand-binding (RLB) models are modified versions of the  Schnakenberg model, where the variations are inspired by biological phenomenon \cite{lee2010aberrant}. The two models are distinguishable from each other in the sense that, in the RLB model, the time delay is incorporated to describe the behavior of both the activator and inhibitor in the system, on the other hand, the LI model includes time delay purely in the dynamics of the activator.

Numerous other studies have investigated time delay in reaction-diffusion systems \cite{ruan1999persistence, ruan1998turing, gourley1996predator, freedman1997global, higham1995existence, chen2013global, chen2011note, chen2013effect, chen2013time, hadeler2007interaction, morita1984destabilization}. Despite this, theoretical analysis is scarce on the basis of the effect of delays on the onset of Schnakenberg--Turing patterns. A thorough investigation of the Schnakenberg model with delay can be found in \cite{yi2017bifurcation}, which confirms the finding that delay has a convoluted relationship with structural robustness. The authors have demonstrated that Hopf bifurcation at the spatially inhomogeneous steady state as a consequence of time delay can induce  Turing instability.

In this paper, we consider a modified form of the Schnakenberg model commonly known as the LI model to investigate complex interactions due to gene expression time delay and domain size in a reaction-diffusion system. In Section \ref{section2}, we introduce the two-dimensional, two-component LI model. In Section \ref{stability_analysis}, we investigate the linear stability analysis to have a theoretical understanding of the effect of time delays and domain length in Turing pattern formation. Section \ref{numerical_analysis} consists of a numerical investigation to explore the dynamics in Turing spaces subject to varying time delays, with the aim to exemplify the results and corroborate the findings of the instability analysis. Additionally, we examine the changes in the Turing instability region due to a change in the domain sizes. Section \ref{pattern_formation} summarizes the findings regarding the initiation of pattern formation within the fixed delay 2D LI model, with a specific emphasis on the influences of delay, domain size, and other parameter values. The paper concludes with a brief summary of our results and discussion in Section \ref{discussion}.

\section{The ligand internalisation (LI) model}
\label{section2}

Ordinary and partial differential equations are powerful tools for modeling complex biological, ecological, physical, and chemical systems. However, many such systems exhibit temporal delays, making models that incorporate their past states more realistic and informative. In this work we consider the 2D LI model, which is a modified version of the Schnakenberg model in the presence of a time delay that represents ligand internalization. The one space-dimensional version was studied in \cite{sargood2022fixed}. The two main rate equations for the activator and the inhibitor are 

\begin{align}
	\begin{aligned}
		\label{eq:Ld}
		\frac{\partial u}{\partial t}&=\frac{d_{u}}{L_x^2}\frac{\partial^2u}{\partial x^2} + \frac{d_{u}}{L_y^2}\frac{\partial^2u}{\partial y^2} + a-u-2u^2v+3\hat{u}^2\hat{v}, \\
		\frac{\partial v}{\partial t}&=\frac{d_{v}}{L_x^2}\frac{\partial^2v}{\partial x^2} + \frac{d_{v}}{L_y^2}\frac{\partial^2v}{\partial y^2} + b-u^2v,\\   
	\end{aligned}
\end{align} 
where $\hat{u} = u(x,y, t-\tau)$, and $\hat{v} = v(x,y, t-\tau)$ signifies activator and inhibitor terms where delay has been considered,  $\tau$ stands for the time delay introduced in the system kinetics, $d_{u}$ and $d_{v}$ represent the diffusion coefficients of the activator $u$ and inhibitor $v$,  respectively. The parameters $d_{u}$ and $d_{v}$ are consistently assigned values of $0.01$ and $0.2$, respectively, throughout the calculation, unless specified otherwise. Parameters $L_x$ and $L_y$ represent the domain size of the rectangle, and after rescaling $[0,L_x]\times[0,L_y]$ into $[0,1]^2$, they become non-dimensional scaling factors of the domain length which we vary to study the outcome of the final pattern, such as stripes or circles. Hence, we consider System \eqref{eq:Ld} on a unit square for our calculations, and varying $L_x$ and $L_y$ will account for the dynamics on different domain sizes. Assuming that morphogens cannot escape at the domain boundaries, Neumann boundary conditions are the most suitable for a biological environment. Thus, 

\begin{align}
\begin{aligned}
\label{bc}
	&\partial_{x} u(0,y,t) = \partial_{x} u(1,y,t) = 
	\partial_{y} u(x,0,t) = \partial_{y} u(x,1,t) = 0,\\
	&\partial_{x} v(0,y,t) = \partial_{x} v(1,y,t) = 
	\partial_{y} v(x,0,t) = \partial_{y} v(x,1,t) = 0.
 \end{aligned}
\end{align}
System \eqref{eq:Ld} has a spatially homogeneous equilibrium solution $\mathcal{E_*} = (u_*,v_*)$, where 

\begin{align}
	\begin{aligned}
		\label{steadystate}
		u_* = a + b, ~~~~ v_* = \frac{b}{(a+b)^2}. 
	\end{aligned}   
\end{align}

% \quad x,y \in\Omega , ~~ t \geq 0

Understanding how delays influence symmetry breaking from a homogeneous steady state is essential, particularly in their interaction with kinetics and diffusion during self-organization. Using analytical and numerical methods, this work evaluates the impact of delays on pattern formation and emphasizes the critical role of domain size. The final heterogeneous steady state of a diffusion-driven instability results from the interplay between reaction, diffusion, and domain geometry, as in  \cite{barrio1999two}, but patterning is also highly influenced by the time delay.

\section{Stability Analysis for Inhomogeneous LI Model}
\label{stability_analysis}
\noindent

In this section, we conduct a detailed analysis of the model \eqref{eq:Ld}-\eqref{bc}, regarding the stability of the spatially uniform steady state, denoted by $(u_*, v_*)$ in \eqref{steadystate}. At the steady state, $u=\hat u=u_*, v=\hat v=v_*$, and the system satisfies the conditions 

\begin{align}
	\begin{aligned}
	f(u_*, v_*,u_*, v_*) = 0,\\
	g(u_*, v_*,u_*, v_*) = 0,
	\end{aligned}
	\label{eq:LI1}
\end{align}
where $f(u, v, \hat u, \hat v)$ and $g(u, v, \hat u, \hat v)$ denote the reaction part of \eqref{eq:Ld} for $u$ and $v$ respectively. To assess the linear stability of the system, we consider small perturbations, $u(x, y, t) = u_* + \delta \zeta(x, y, t)$, $v(x, y, t) = v_* + \delta \eta(x, y, t)$, where $|\delta| \ll 1$, around the steady state. 
Expanding $f$ and $g$ in the Taylor series and truncating at $\mathcal{O}(\delta)$, we find that perturbations evolve according to the linearized system

\begin{align}
	\begin{aligned}
		\begin{cases}	
		&\frac{\partial \zeta}{dt} = \frac{d_{u}}{L_x^2} \frac{\partial^2\zeta}{\partial x^2} + \frac{d_{u}}{L_y^2} \frac{\partial ^2\zeta}{\partial y^2} - \zeta - 4 u_*v_* \zeta - 2u_*^2\eta+ 6u_*^2v_*\hat\zeta+3u_*^2\hat\eta ,\\
		&\frac{\partial \eta}{dt}=  \frac{d_{v}}{L_x^2} \frac{\partial^2\eta}{\partial x^2} + \frac{d_{v}}{L_y^2} \frac{\partial^2\eta}{\partial y^2}- 2u_*v_*\zeta -u_*^2\eta,\\   &\zeta_x(0,y,t)=\zeta_x(1,y,t)=\zeta_y(x,0,t)=\zeta_y(x,1,t)=0,\\
		&\eta_x(0,y,t)=\eta_x(1,y,t)=\eta_y(x,0,t)=\eta_y(x,1,t)=0,\\
		\end{cases}
	\end{aligned}
	\label{eq:linear}
\end{align}
where $\hat\zeta=\zeta(x,y,t-\tau)$, $\hat\eta=\eta(x,y,t-\tau)$. We now look for solutions of the form\\

$$(\zeta(x,y,t),\eta(x,y,t))^\top = e^{\lambda t}\cos(k_x\pi x)\cos(k_y\pi y)(c, d)^\top, $$ 

where $k_x,k_y\in\mathbb{N}_0$, $c,d\in \mathbb{R}$. Substituting the above into \eqref{eq:linear}, we have:

\begin{align}
	\begin{aligned}
		\lambda c =& -\frac{d_{u} k_x^2 \pi^2}{L_x^2} c  -\frac{d_{u} k_y^2 \pi^2}{L_y^2} c - c - 4u_* v_* c + 6 u_*v_* e^{-\tau \lambda} c -2u_*^2 d + 3 u_*^2e^{-\tau \lambda d},\\
		\lambda d =& - 2u_* v_*c  -\frac{d_{v} k_x^2 \pi^2}{L_x^2}  d -\frac{d_{v} k_y^2 \pi^2}{L_y^2} d - u_*^2 d.\\	
	\end{aligned}
	\label{eq:linear1}
\end{align}

The above system has a nonzero solution for $(c, d)$ if and only if

\begin{equation}
	\det \begin{bmatrix}
		-\frac{d_{u} k_x^2 \pi^2}{L_x^2}  -\frac{d_{u} k_y^2 \pi^2}{L_y^2} -1 - 4u_* v_* + 6 u_*v_* e^{-\tau \lambda} - \lambda   &-2u_*^2 + 3 u_*^2e^{-\tau \lambda}   \\
		- 2u_* v_*&  -\frac{d_{v} k_x^2 \pi^2}{L_x^2}  -\frac{d_{v} k_y^2 \pi^2}{L_y^2} - u_*^2 - \lambda \\
	\end{bmatrix} = 0_{\textstyle\raisebox{3pt}{.}}
\end{equation}

Hence, we obtain the characteristic equation,

\begin{align}
	\begin{aligned}
		\label{eq:ch1}
		D_{k_{x}k_{y}}(\lambda,\tau)&= \lambda_k^2 + p_{k_x,k_y} \lambda + q_{k_x,k_y} + (r_{k_x,k_y}\lambda + s_{k_x,k_y})e^{-\lambda \tau} = 0, ~ k_x,k_y \in \mathbb{N}_0,
	\end{aligned}
\end{align}	

where

\begin{align}
	\begin{aligned}
		p_{k_x,k_y} =&\frac{(d_{u} + d_{v}) k_{x}^2 \pi^2}{L_{x}^2} + \frac{(d_{u} + d_{v}) k_{y}^2 \pi^2}{L_{y}^2} + u_*^2 + 4 u_* v_* + 1,\\
		q_{k_x,k_y} =& d_{v} \left (\frac{k_{x}^2}{L_{x}^2} + \frac{k_{y}^2}{L_{y}^2} \right) \pi^2 + d_{u} d_{v} \left (\frac{k_{x}^4}{L_{x}^4} + \frac{k_{y}^4}{L_{y}^4} \right) \pi^4 + d_{u} \left (\frac{k_{x}^2}{L_{x}^2} + \frac{k_{y}^2}{L_{y}^2} \right) \pi^2 u_*^2\\& + 4 d_{v} \left (\frac{k_{x}^2}{L_{x}^2} + \frac{k_{y}^2}{L_{y}^2} \right) \pi^2 u_* v_* + \frac{2 d_{u} d_{v} k_{x}^2 k_{y}^2 \pi^4}{L_{x}^2 L_{y}^2} + u_*^2,\\
		r_{k_x,k_y} =& -6u_*v_*,\\
		s_{k_x,k_y} =& -\frac{6 d_{v} k_{x}^2 \pi^2 u_* v_*}{L_{x}^2} - \frac{6 d_{v} k_{y}^2 \pi^2 u_* v_*}{L_{y}^2}.\\ 
	\end{aligned}
\end{align}	

This dispersion relation \eqref{eq:ch1} provides the function of eigenvalues, represented by $\lambda$, corresponding to different wave numbers, denoted by $k_{x}$ and $k_{y}$. 
The parameter set $(a, b, d_{u}$, $d_{v}, \tau, L_{x}, L_{y})$ responsible for the existence of solutions in order for patterns to arise can be computed from the characteristic equation \eqref{eq:ch1}. The term $e^{-\lambda \tau}$ captures the primary influence of delay on linear stability. We introduce the following notation for the maximal real parts of the eigenvalues for different values of $k_x$ and $k_y$ as 

\begin{align}
	\begin{aligned}
		\label{eq:alpha1}
	\alpha &= \max\{\alpha_{k_x, k_y}(\tau) : \ k_x, k_y \in \mathbb{N}_0\},
	\end{aligned}
\end{align}

where

\begin{align}
	\begin{aligned}
		\label{eq:alpha2}
		\alpha_{k_x, k_y} (\tau) & = \max\{\mathrm{Re}\lambda:\  D_{k_{x}k_{y}}(\lambda,\tau) = 0 \}, \quad k_x, k_y \in \mathbb{N}_0.
	\end{aligned}
\end{align}

The linear theory suggests that a general perturbation will grow with rate  $ e^{ \alpha(\tau) t}$, given that $\alpha(\tau) > 0$, and hence the emergence of patterns can be anticipated. This theory can also be utilized to study the implications that the introduction of time delays and varying the domain size into the reaction kinetics can possibly have. To accomplish this, we illustrate bifurcation diagrams (otherwise referred to as Turing space or $(a,b)$ parameter plane) to signify the regions corresponding to Turing instability where pattern formation emerges. The Turing space allows us to choose parameters $a,b$ which facilitates the construction of plots depicting the relationship between $\alpha$ and $\tau$, as well as plots illustrating $\alpha_{k_x, k_y}$ with respect to variations in $L_{x}$. These visualizations offer an initial insight into the behavior of the reaction-diffusion system.

\subsection{The $(a,b)$ parameter plane}

\begin{figure}[H]
	\centering
    \includegraphics[width=0.4\linewidth]{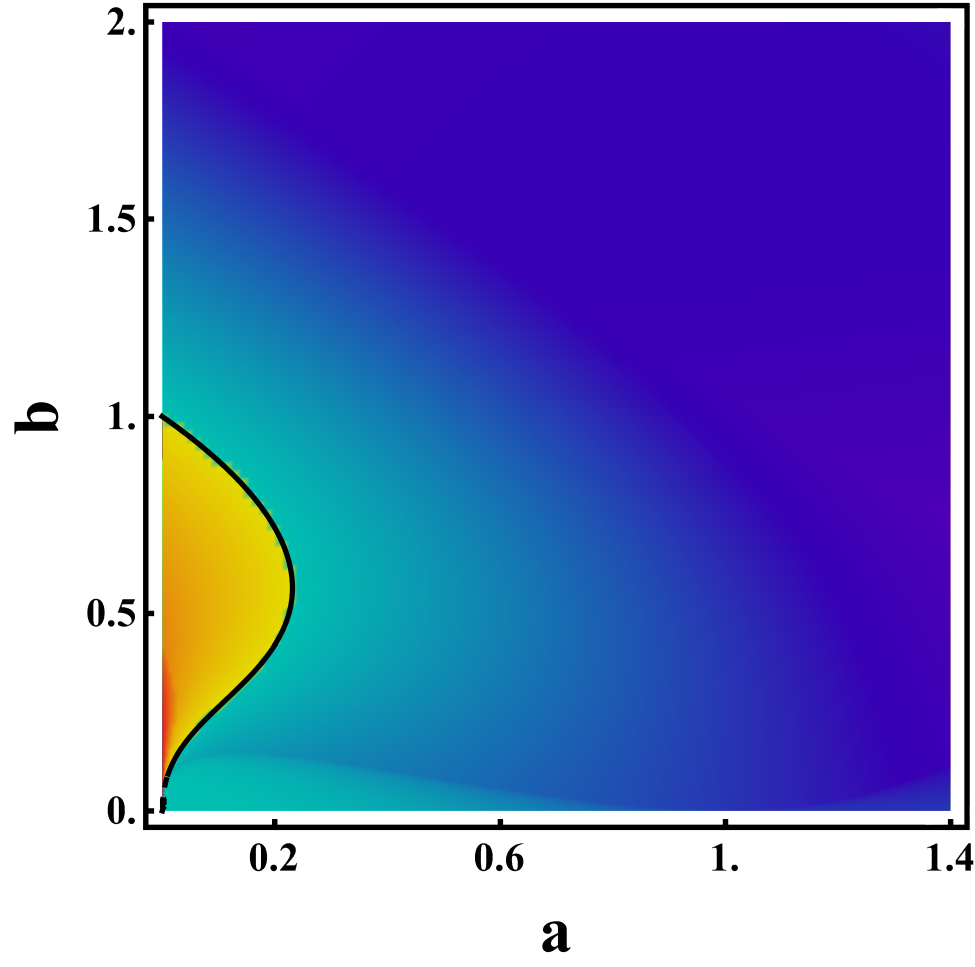}
    \includegraphics[width=0.4\linewidth]{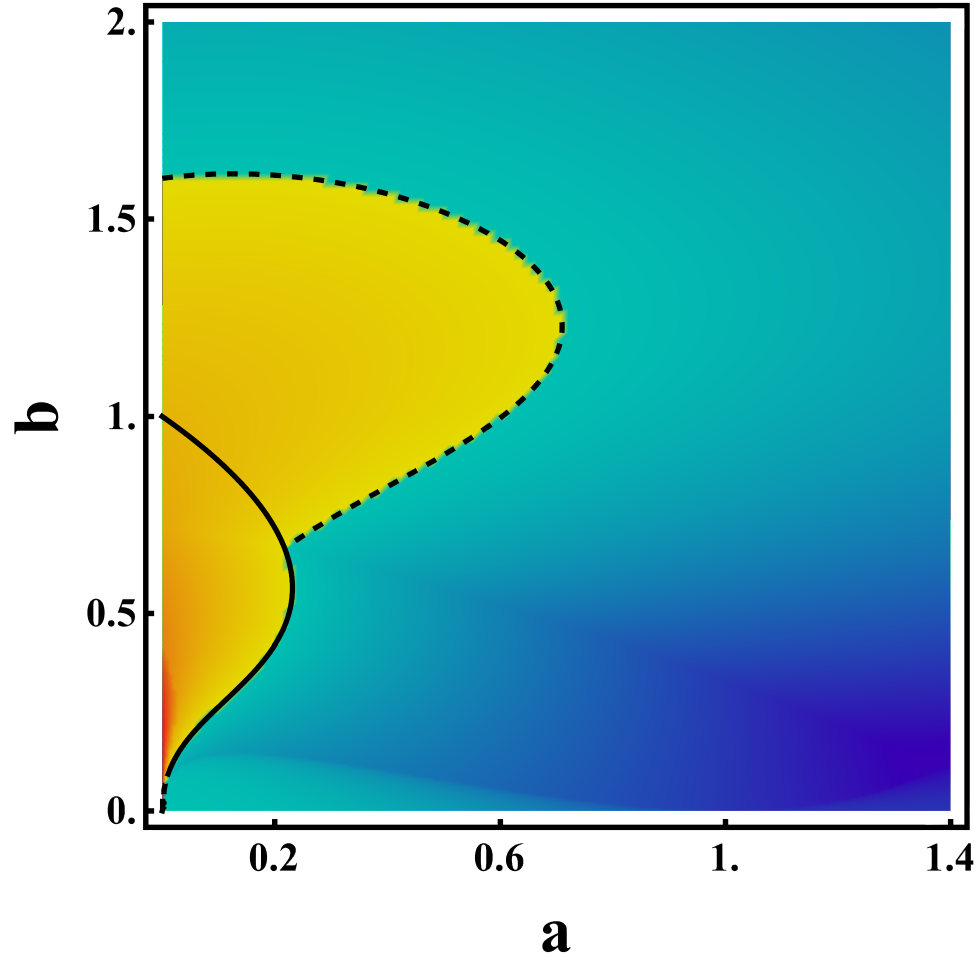}
	\raisebox{0.4 cm}{\includegraphics[height=6cm]{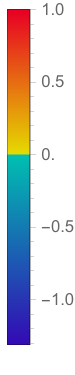}}\\
    \includegraphics[width=0.4\linewidth]{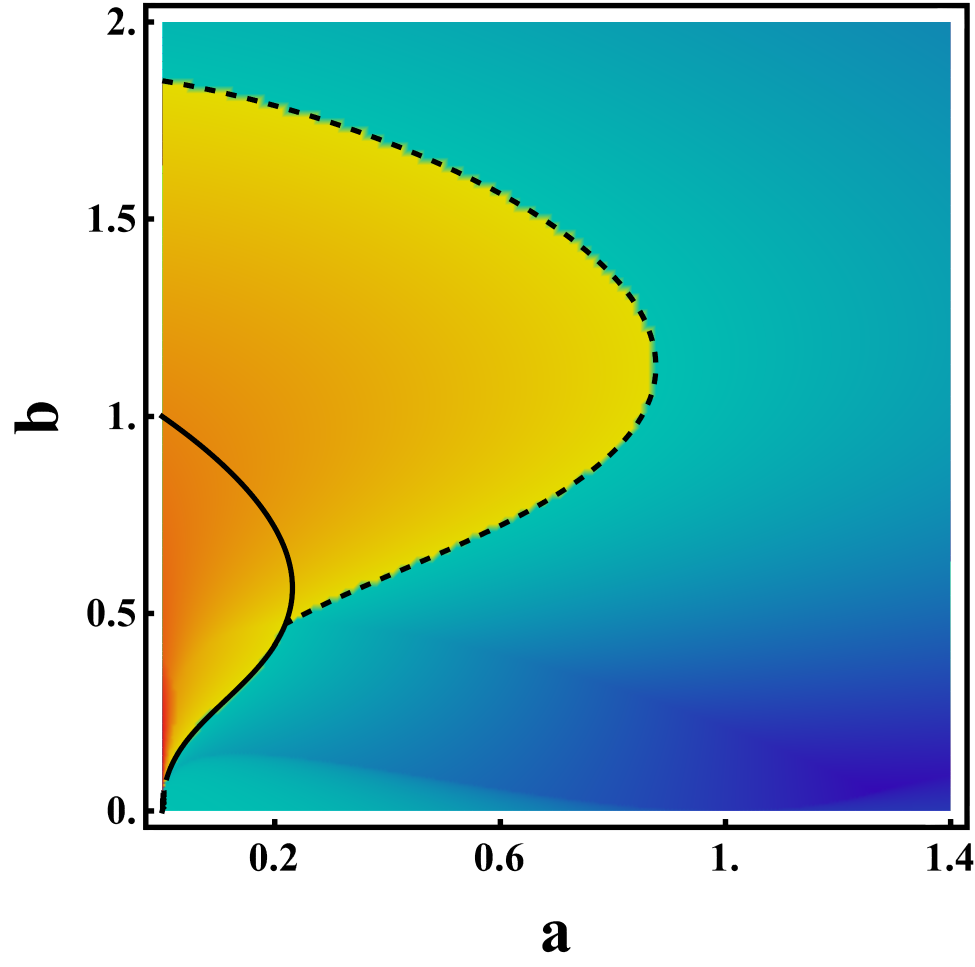} 
    \includegraphics[width=0.4\linewidth]{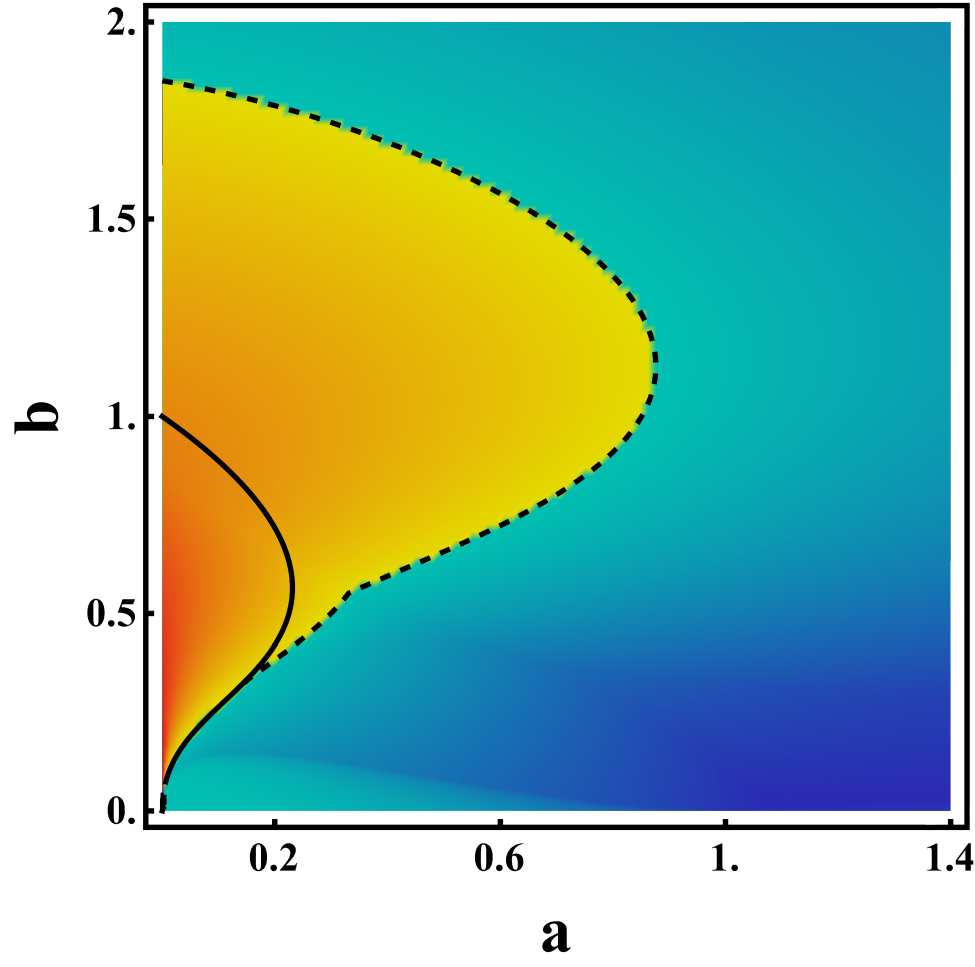}
	\raisebox{0.4 cm}{\includegraphics[height=6cm]{fig_3_scale.pdf}}
	\caption{Heatmap of the maximal real part of the eigenvalues, $\alpha$, as defined in \eqref{eq:alpha1}, computed over $a$ and $b$ parameter space by solving the characteristic equation (\ref{eq:ch1}) without delay ($\tau = 0$). The stability curves for the spatially homogeneous (solid black curve) and spatially inhomogeneous (black dotted curve) conditions are shown. The area between these two curves signifies the Turing region. Here $d_{u} = 0.01$, $d_{v} = 0.2$ and $L_{y} = 0.2$. For (a)--(d), $L_{x} = 0.2, 0.38, 0.5, 1$ respectively. As $L_{x}$ increases, the Turing space expands.}
	\label{fig:fig3}
\end{figure}

In \cref{fig:fig3} (a--d), we present stability diagrams depicting the behavior of the spatially homogeneous fixed point. These diagrams display the heatmap of $\alpha$, as defined in equation \eqref{eq:alpha1}, across the parameter space defined by $(a, b)$ in the absence of delay, i.e. $\tau = 0$. We systematically vary the domain size, as depicted in \cref{fig:fig3} (a--d). To accentuate the Turing instability region, we superimpose contour lines corresponding to $\alpha = 0$ and $\alpha_{0,0} = 0$, where $\alpha$ and $\alpha_{0,0}$ are defined by equations \eqref{eq:alpha1} and \eqref{eq:alpha2}, respectively.

In \cref{fig:fig3} (a), where $L_x = 0.2$, the spatial extent is insufficient to facilitate the formation of the Turing instability region, therefore the curve $\alpha = 0$ is absent. Conversely, in \cref{fig:fig3} (b), (c), and (d), as the domain length $L_{x}$ is gradually increased to values of $0.38, 0.5$, and $1$ respectively, two distinct curves become evident. The black dotted curve delineates the boundary of instability for the fixed point in the presence of diffusion, and, the solid black curve outlines the boundary for the region where the fixed point remains unstable, even in the absence of diffusion. The region enclosed between these two curves can be defined as the "Turing space." This designation signifies the parameter space wherein Turing patterns may potentially form. 

The figures provide a visual representation of the conditions where the fixed point's behavior changes, under the effects of varying domain sizes, with the two different curves representing the stability or instability of a fixed point under the influence of diffusion and in its absence. With an increase in domain size, there is a distinct expansion of the region associated with Turing instability, illustrated by the black dotted curve associated with $\alpha = 0$. This expansion widens the available space where patterns may emerge. We select parameter values for $a$ and $b$ from this graphical representation and with the advantage of numerical simulations using these values, we get further insight into the dynamic behavior associated with Turing patterns. It is also worth mentioning, an inflection point is visible in the black dotted curve when $L_{x}$ attains a value of 1, which could imply a change in the dominant mode. Note that the parameter $L_{y}$ is consistently fixed at $0.2$ to generate this figure. It is also essential to highlight that $d_{u}$ and $d_{v}$ are typically set to $0.01$ and $0.2$ respectively unless explicitly stated otherwise.

\begin{figure}[H]
	\centering
	\includegraphics[width=0.4\linewidth]{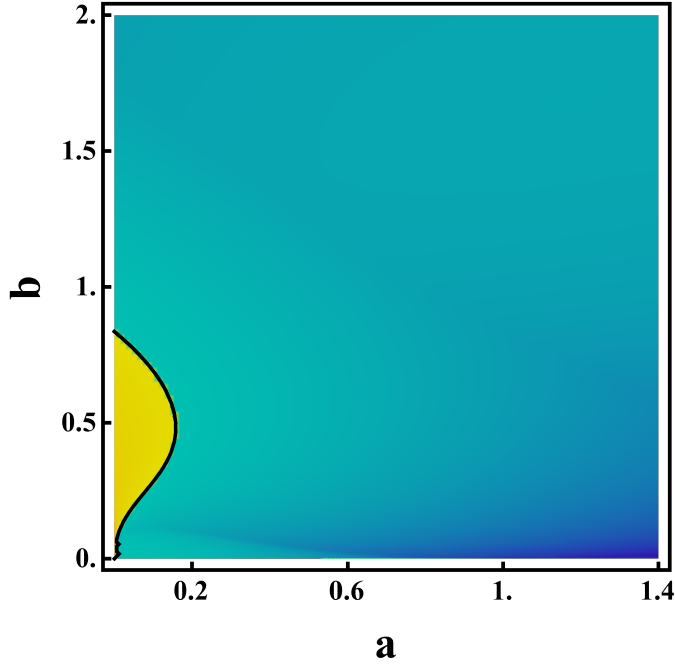}
    \includegraphics[width=0.4\linewidth]{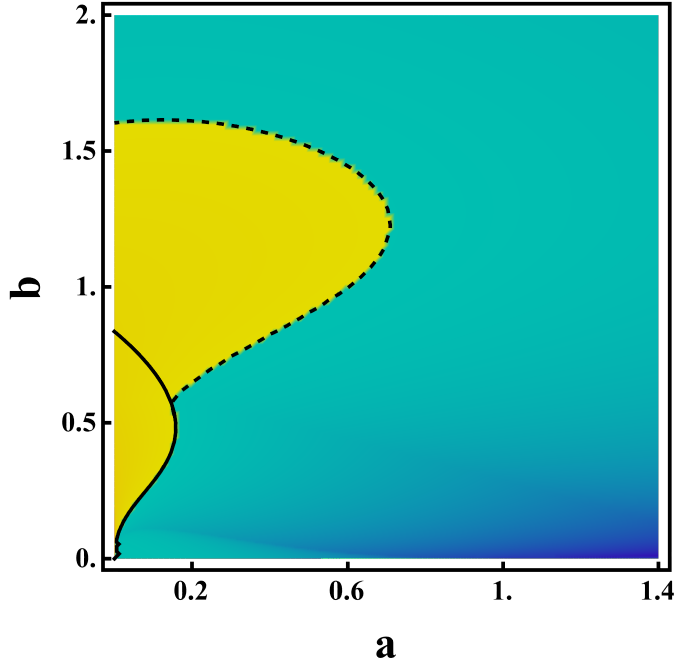}
	\raisebox{0.4 cm}{\includegraphics[height=6cm]{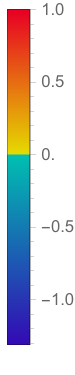}}\\
    \includegraphics[width=0.4\linewidth]{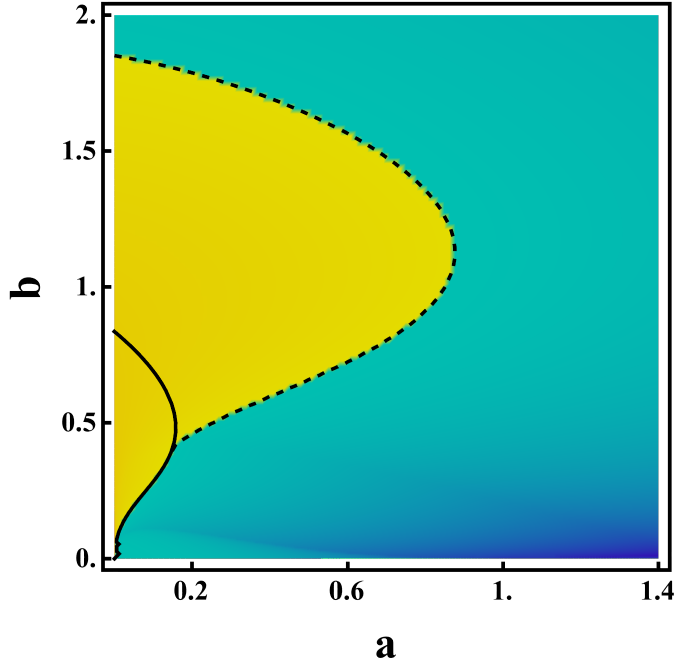} 
    \includegraphics[width=0.4\linewidth]{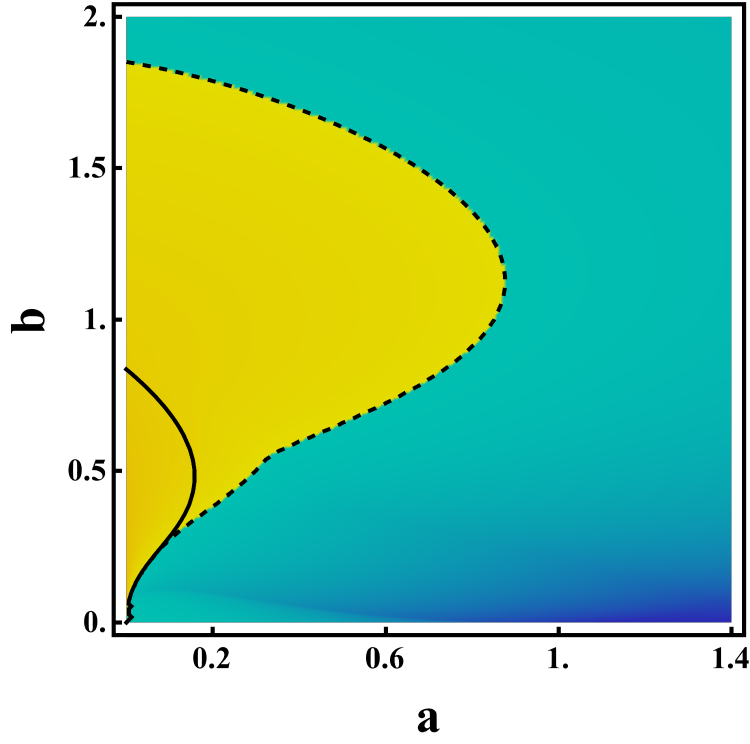}
	\raisebox{0.4 cm}{\includegraphics[height=6cm]{fig_4_scale.pdf}}
	\caption{(a)--(d) $\alpha_{k_x,0}$ as defined in \eqref{eq:alpha1} is plotted as a function of $L_{x}$ for different modes for equation \ref{eq:ch1} for a fixed delay. In (a), $L_{y} = 0.2$ and (b), $L_{y} = 3$; $\tau = 0$ in both cases. In (c), $L_{y} = 0.2$ and (d) $L_{y} = 3$; $\tau = 1$ in both cases. $a = 0.1, b = 0.9$, $d_{u} = 0.01$, and $d_{v} = 0.2$ through (a--d). A significant variability in the dominant mode is noticed for smaller domain sizes. As we increase $\tau$ to 1, the eigenvalues tend toward zero, leading to a higher density of lines.}
	\label{fig:fig2}
\end{figure}

In \cref{fig:fig2}, we present a similar stability diagram, but with a constant delay of $\tau = 1$ in the system. From the figure, it is evident that as we increment the time delay $\tau$ to 1, the region of instability in the absence of diffusion (bounded by the solid curve $\alpha_{0,0} = 0$), decreases in size, potentially enlarging the Turing instability region. As $\tau$ increases the absolute magnitude of $|\alpha|$ decreases.  This suggests that the eigenvalues draw closer to the origin and hence, the emergence of patterns will necessitate more time when larger delays are introduced. Enlarging the domain, by comparing (a)-(d), we observe that increasing domain size extends the Turing instability region.

\subsection{Plot of $\alpha$ against $\tau$}

\begin{figure}[]
	\begin{center}
		\includegraphics[width=0.8 \linewidth]{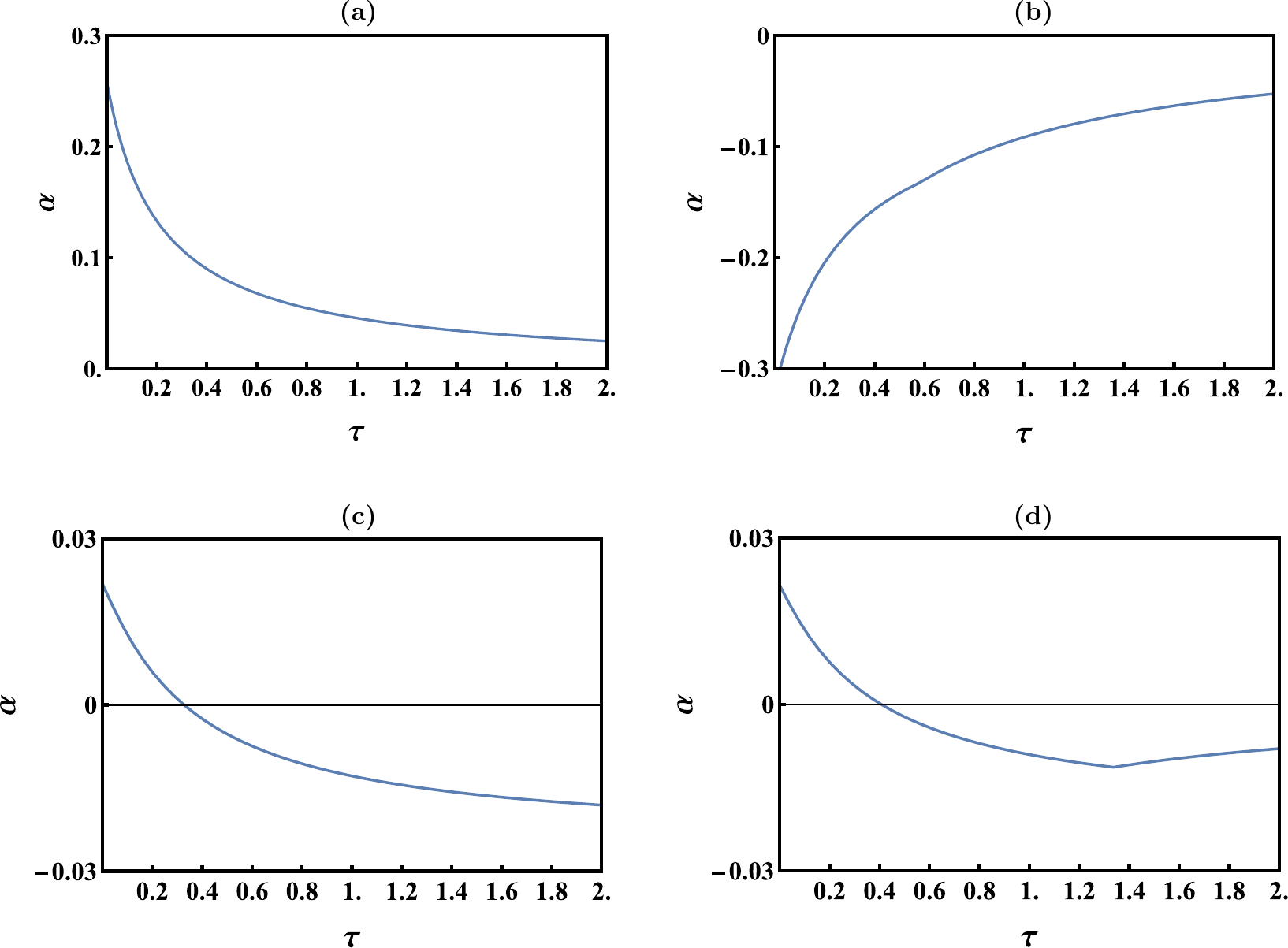}\\
		\caption{Plot of the maximal real part of the eigenvalues against time delay. In (a)--(d), $\alpha$ as defined in \eqref{eq:alpha1} is evaluated for fixed delay for \eqref{eq:ch1} as $\tau$ varies from $0$ to $2$. In the case (a), $(a, b) = (0.1, 0.9)$, pattern formation can be predicted (since $\alpha > 0$). Conversely, in case (b), $(a, b) = (0.4, 0.4)$, no pattern formation can be anticipated (since $\alpha < 0$). In both (a) and (b), $L_{y} = 0.2$, $L_{x} = 3$. In (c), $(a, b) = (0.17, 0.5)$, $L_{y} = 0.2$, $L_{x} = 0.38$. In (d), $(a, b) = (0.18, 0.4)$, $L_{y} = 0.2$, $L_{x} = 0.5$. In both (c) and (d), equilibrium becomes stable after a certain $\tau$. In (d), a mode switch is also apparent at $\tau \approx 1.33$. $d_{u} = 0.01$, and $d_{v} = 0.2$ is fixed through (a)--(d).}\label{fig:fig1}
	\end{center}
\end{figure}

In this section, we generate plots for real parts of the eigenvalues to investigate possible occurrences of Turing instability for $(a, b)$ parameter pairs chosen from stability diagrams, \cref{fig:fig3} and \cref{fig:fig2}. In \cref{fig:fig1} (a)--(d), the blue curve corresponds to different values of 
$\alpha$, as defined in \eqref{eq:alpha1} plotted against $\tau$ for four different combinations of parameters $a$ and $b$. \cref{fig:fig1}(a), shows an example where Turing instability occurs when the parameter values are set to  $(a,b) = (0.1, 0.9)$. As $\tau$ increases, $\alpha$ decreases but remains positive. This indicates that pattern formation will consistently occur with this particular combination of $a$ and $b$. On the other hand, in \cref{fig:fig1}(b), with parameter values set to  $(a,b) = (0.4, 0.4)$, as $\tau$ increases, the values of $\alpha$ grow larger, yet they remain negative. This suggests that, for the chosen parameter set, pattern formation does not occur, and the spatially homogeneous equilibrium remains stable.

In \cref{fig:fig1} (c) and (d), we explore scenarios for pattern formation for other combinations of $a$ and $b$ values. This analysis provides clarification of the intricate relationship between pattern formation and domain size. In Subfigure \ref{fig:fig1} (c), a noticeable trend emerge as $\tau$ takes on smaller values, particularly for the specific parameter combination of $(a, b) = (0.17, 0.5)$, within the confines of relatively modest domain sizes, denoted by $L_{x} = 0.2$ and $L_{y} = 0.38$. As $\tau$ exceeds the critical threshold of approximately $0.32$, the manifestation of pattern formation becomes unidentifiable. \cref{fig:fig1}(d) shows a parallel outcome when considering the parameter values $(a, b) = (0.18, 0.4)$, in conjunction with domain dimensions $L_{x} = 0.2$ and $L_{y} = 0.5$. Significantly, as $\tau$ approaches a value of approximately $0.39$, the system achieves a state of stability. A switch in the dominant mode is noted; this transition however does not yield any noticeable impact on the overall system dynamics.

In overview, \cref{fig:fig1} provides valuable insights into the significant impact that system parameters, particularly $(a, b, \tau, L_{x}, L_{y})$, can exert on the system's dynamics. It provides an initial understanding of how the values of $(a, b)$ can influence the system and the role that the system's history and spatial dimensions, represented by $L_{x}$ and $L_{y}$, play in shaping the dynamic changes. Numerical simulations in later sections are conducted based on \cref{fig:fig1} (a).

\subsection{Plot of $\alpha_{k_x, k_y}$ against $L_{x}$}

\begin{figure}[H]
	\centering
	\includegraphics[width=0.4\linewidth]{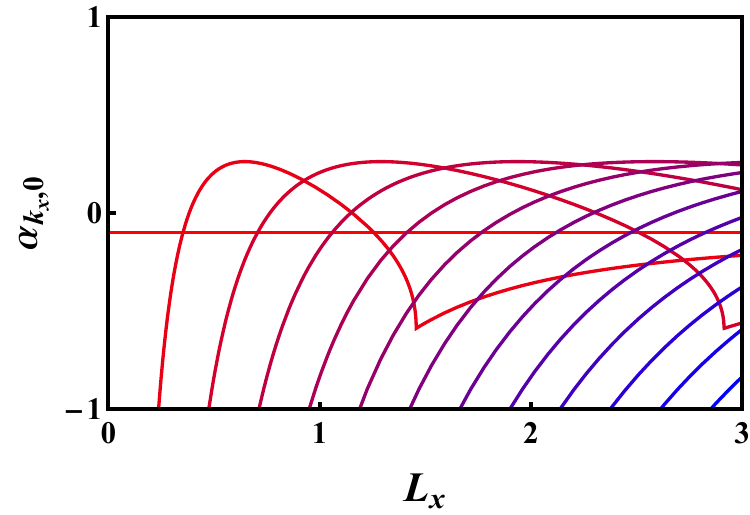}  
	\includegraphics[width=0.4\linewidth]{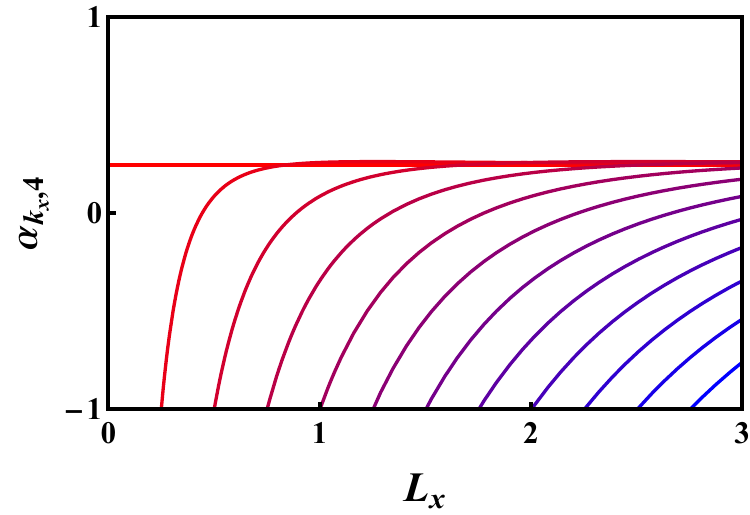}
	\includegraphics[width=0.4\linewidth]{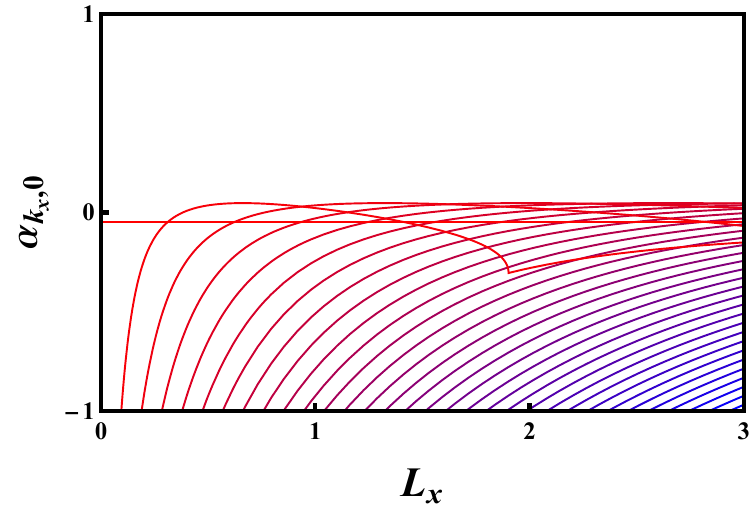}
	\includegraphics[width=0.4\linewidth]{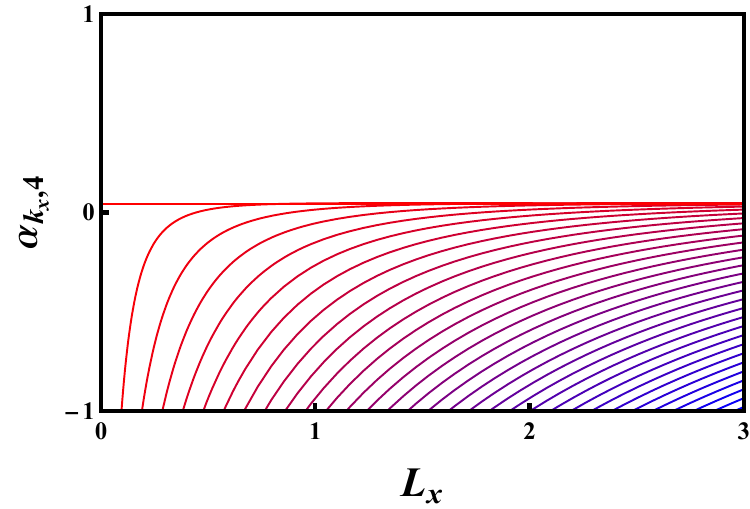}
	\caption{(a)--(d) $\alpha_{k_x,0}$ as defined in \eqref{eq:alpha1} is plotted as a function of $L_{x}$ for different modes for equation \ref{eq:ch1} for a fixed delay. In (a), $L_{y} = 0.2$ and (b), $L_{y} = 3$; $\tau = 0$ in both cases. In (c), $L_{y} = 0.2$ and (d) $L_{y} = 3$; $\tau = 1$ in both cases. $a = 0.1, b = 0.9$, $d_{u} = 0.01$, and $d_{v} = 0.2$ through (a--d). A significant variability in the dominant mode is noticed for smaller domain sizes. As we increase $\tau$ to 1, the eigenvalues tend toward zero, leading to a higher density of lines.}
	\label{fig:fig4}
\end{figure}

In \cref{fig:fig4}, we generate plots for $\alpha_{k_x, k_y}$ defined by \eqref{eq:alpha2} as functions of $L_x$. These graphs are relevant for analyzing the behavior of the dominant eigenvalue with respect to variations in $L_{x}$, the domain size. In \cref{fig:fig4} (a), the choice of a relatively small value for $L_{y}$, specifically set at $0.2$, leads to the dominance of the mode with an index of $k_y = 0$. The horizontal line corresponds to the scenario where both $k_x$ and $k_y$ are equal to zero. In this case, as the parameter $L_x$ undergoes variation, the value of $\alpha_{0,0}(L_x)$ remains constant and consistently negative. This observation implies that in the absence of diffusion, the spatially homogeneous equilibrium remains stable. Conversely, when examining $\alpha_{k_x,0}(L_x)$ where $k_x > 0$, as we increment $L_x$, $\alpha_{k_x,0}(L_x)$ exhibits a continuous upward trend until it attains its maximum positive value. Subsequently, it sustains a continuous downward trend until it reaches its minimum negative value. At this point, two real eigenvalues collide and transform into a pair of complex conjugate eigenvalues. It then follows an upward trajectory once more until it ultimately converges at a negative value. As a result, we can assert that the predominant mode modifies as a consequence of our deliberate selection of a reduced domain size in the context of $L_{y}$. In \cref{fig:fig4} (b), with $L_{y} = 3$, it is evident that a shift in the dominant mode is not a prevalent occurrence when the domain size is substantially large. Note that, in \cref{fig:fig4} (a) and (b), the delay parameter $\tau$ is set to $0$. 

We introduce a comparable set of results, now considering the presence of a delay in the system, with $\tau = 1$ in \cref{fig:fig4} (c) and (d). As we increment the value of $\tau$, the eigenvalues shift closer to the origin. This phenomenon explains the occurrence of increased density of lines in the figures. As before, it is evident that when the domain length, denoted as $L_{y}$, is relatively small, there is a distinct and pronounced transition in the dominant mode. However, when we extend the domain length to $L_{y}=3$, the transition becomes imperceptible. 

The above observations underscore the sensitivity of the system to variations in domain length, with shorter domains exhibiting a more pronounced mode switch, while longer domains yield a less noticeable transition. Furthermore, it is apparent that when the parameter $\tau$ is subjected to increment, the eigenvalues converge towards the origin, indicative of a distinct alteration in the system's dynamics.  Supplementary results are provided in the subsequent sections to substantiate these phenomena.

\section{Simulation}
\label{numerical_analysis}

In this section, we focus on the numerical analysis of the dynamic behavior of the system presented in \Cref{eq:Ld}, which is a delayed form of the Schnakenberg model. Understanding the delay is of paramount importance in biological systems, as is the inclusion of delays in mathematical models incorporating biological phenomena. For example, the delay in protein processing of Nodal-related proteins, such as Cyc and Sqt, which are instrumental during various developmental stages in zebrafish, has been investigated in \cite{jing2006mechanisms}. The interconnection of a biological oscillator, with focus on the critical role of delay mechanisms in shaping system behavior was studied in \cite{lewis2003autoinhibition}. By investigating the somitogenesis oscillator in zebrafish embryos, the paper delves into how delays in gene expression processes influence the timing and synchronization of oscillatory patterns. 

The role of delay in biological processes is complex and remains an active area of research. Using our model, we analyze how gene expression time delays affect system dynamics and pattern formation, with particular attention to domain size. Our results reveal that small domain sizes can strongly influence pattern dynamics, and the delay and domain size together have an intricate effect on patterning.

We carry out computational simulations to investigate the dynamics of pattern formation in 2D Schnakenberg kinetics in the presence of delay. With the help of the results obtained in the linear instability analysis, we also investigate and compare the relationship between the time for the pattern to arise with $\tau$ and $L_{x}$. To perform the simulations, we consider $n \times m$ grid and pay special attention to a system where the length is rescaled to 1, i.e. $\Delta x = 1/(n-1)$ and $\Delta y = 1/(m-1)$ ensuring $x_1, y_1 = 0$ and $x_n, y_m = 1$.
\subsection{Initial conditions for numerical simulation}
\label{section4}
To solve system \eqref{eq:Ld}, $u = u(x,y, t)$, and $v = v(x,y, t)$, the concentrations of the activator and inhibitor respectively at the spatial location $(x,y)$ are subjected to initial conditions:
\begin{align}
	\begin{aligned}
		u(x, y, 0) &= u_* + \mathcal{R}_1(x, y), \\
		v(x, y, 0) &= v_* + \mathcal{R}_2(x, y),\\
		t \in [-\tau, 0],\ x &\in [0, L_{x}],\ y \in [0, L_{y}],
	\end{aligned}\label{initial_1}
\end{align} 
where $u_* = a + b$ and $v_* = \frac{b}{(a+b)^2}$ are the steady-state solutions of the system, $\mathcal{R}_1$ and $\mathcal{R}_2$ stand for small random perturbations.

The spatial derivatives were discretized using the method of finite difference on an $n \times m$ grid, using the conventional five-point stencil for the Laplacian to evaluate the second-order derivative in the $x$ and $y$ direction, namely,

\begin{align}
	\begin{aligned}
		\frac{\partial^2u}{\partial x^2} (x,t) &= \frac{1}{\Delta x^2} (u_{i+1} - 2u_{i} + u_{i-1}),\\
		\frac{\partial^2u}{\partial y^2} (y,t) &= \frac{1}{\Delta y^2} (u_{j+1} - 2u_{j} + u_{j-1}).
	\end{aligned}
\end{align}

The stability of the numerical simulation is assured by implementing the Courant--Fried\-rich--Lewy condition, which is a necessary convergence condition for numerically solving such partial differential equations. This stability condition states that for information to pass through a space discretization, the time step $\delta t$ should satisfy $\max\{d_u, d_v\}\left(\frac{1}{(\delta x)^2}+\frac{1}{(\delta y)^2}\right) \delta t \le \frac {1}{2}$, where $\delta x$ and $\delta y$ are space steps.

For comparative analysis, we also perform a second set of calculations by considering the initial conditions as a function of $cosine$ shown here below:
\begin{align}
	\begin{aligned}
		u_0(x,y,t) &= u_* + e^{\lambda t} c \cos(k_{x} \pi x) \cos(k_{y} \pi y),\\
		v_0(x,y,t) &= v_* + e^{\lambda t}d \cos(k_{x} \pi x) \cos(k_{y} \pi y),
	\end{aligned}\label{initial_2}
\end{align} 
where $u_*$ and $v_*$ are the steady-state solutions of the system, $(c, d)^\top$ is the eigenvector of the matrix, $k_{x}$ and $k_{y}$ are the indices of the dominant mode. The dominant modes can be computed from the characteristic equation for varying $L_{x}$ and $L_{y}$ values. This initial condition has an unequivocal correlation with the characteristic equation. 

\subsection{Time to Pattern }
One of the primary consequences resulting from the inclusion of delays in the reaction-diffusion equation is its influence on the initiation of pattern formation. In this particular context, the term "time to pattern"  was first mentioned in \cite{sargood2022fixed} to represent the time duration required for patterns to emerge. For the sake of simplicity, we adopt similar nomenclature. It is essential to emphasize that in our case, the measurement of the time to pattern encompasses not only various delays but also different domain sizes. The time to pattern, which we denote by $\mathcal{T}$, is measured as the first time when the perturbation from the steady state is significantly large in the maximum norm, such that it reaches a magnitude beyond a certain threshold, signifying the onset of pattern formation.

To quantify time to pattern, we initiate a function perturbed from the steady state such as, 

\begin{align}
	\begin{aligned}
		u(x,y) - u^{*} &\approx e^{\lambda t} c \cos(k_{x}\pi x)\cos(k_{y}\pi y), \\
		v(x,y) - v^{*} &\approx e^{\lambda t} d \cos(k_{x}\pi x)\cos(k_{y}\pi y).
	\end{aligned}\label{eq:perturb}
\end{align}

Here $u^{*}$ and $v^{*}$ represent the steady-state values of the variables $u$ and $v$ respectively. The distance between the perturbed state $(u_{t},v_{t})$ and the steady state  $(u^{*},v^{*})$ denoted by $||(u_t,v_t) - (u^{*},v^{*})||$ can be approximated as  $e^{\lambda t} ||(u_0,v_0) - (u^{*},v^{*})||$, where $(u_0,v_0)$ is the initial condition. Introducing a threshold in our analysis denoted by $w  =  e^{\lambda t} \beta$, where $\beta = ||(u_0,v_0) - (u^{*},v^{*})||$, we can rearrange the relation to solve for the time $\mathcal{T}$ required for the distance to reach the threshold:

\begin{align}
	\begin{aligned}
		\mathcal{T} \approx \frac{1}{\lambda} \ln \frac{w}{\beta},
	\end{aligned}\label{eq:ttp}
\end{align}
where $\lambda$ corresponds to $\alpha$. To evaluate $\mathcal{T}$ in numerical simulations we consider $\beta = 0.005$ and $w = 0.01$ and also use the same value to determine the $\mathcal{T}$ from the characteristic equation in order to maintain consistency.

\begin{figure}[H]
	\begin{center}
		\includegraphics[width=0.8 \linewidth]{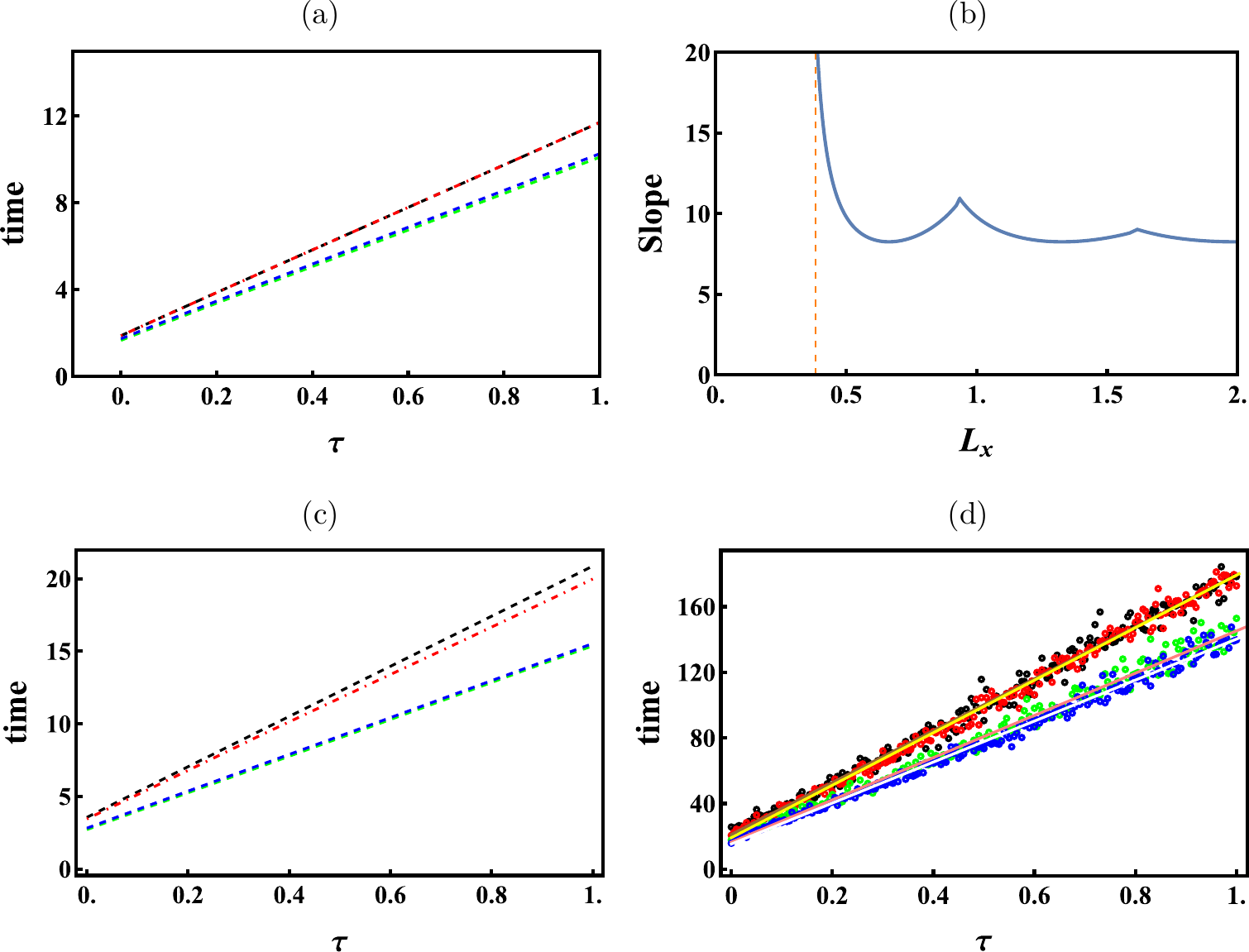}\\
		\caption{(a) Predicted time to pattern with respect to delay computed from the characteristic equation \eqref{eq:ch1} showing a linear trend for different $L_{x}$ = $0.5$ (black, dashed), $1$ (red, dashed), $1.2$ (green, dashed), $1.5$ (blue, dashed). (b) The change in the slope with $L_{x}$ for the curves in (a). (c) Simulated Time to pattern for the initial condition in \eqref{initial_2} for $L_{x} = 0.5$ (black, dashed), $1$ (red, dashed), $1.2$ (green, dashed), and $1.5$ (blue, dashed). (d) Time to pattern simulated considering random initial conditions \eqref{initial_1} for $L_{x} = 0.5$ (black dots), $1$ (red dots), $1.2$ (green dots) and $1.5$ (blue dots). In all the cases, $a=0.1$, $b =0.9$, $d_{u} = 0.01$, $d_{v} = 0.2$, $L_{y}= 0.2$.}\label{fig:fig5}
	\end{center}
\end{figure}

Simulations are subject to altering $\tau$ to explore its impact on the time required for pattern formation. As depicted in \cref{fig:fig5} (a), we found a linear correlation between time to pattern and $\tau$, as evidenced by analysis of the characteristic equation. This result aligns with the earlier-established conclusion by Gaffney et al. \cite{gaffney2006gene}, who demonstrated a linear correlation within the context of a one-dimensional LI model. The linear dependency suggests that the increase in the delay parameter leads to a proportional increase in the time required for the emergence of patterns.

We repeat the evaluation of time to pattern for varying $L_{x}$ while keeping the other parameters constant and compare the solutions. A consistent linear trend is observed in all cases, indicating that the computed value of $\lambda_ {k}$ effectively predicts the time to pattern formation. This observed time variation exhibits a strong linear correlation with $\tau$, demonstrating a high degree of accuracy in approximation. However, it can be seen that resultant slope varies for different value of $L_{x}$. The various curves in \cref{fig:fig5} (a) represent time to pattern corresponding to different domain sizes plotted against the delay $\tau$, and it becomes apparent that the system exhibits a distinct nonlinearity as the domain size varies, as exemplified by the slope presented in \cref{fig:fig5} (b). This outcome highlights intriguing phenomena regarding the impact of domain size on the dynamics of pattern formation, an observation not previously documented in existing studies. As a result, these findings motivated us to conduct investigations into the domain size.

We substantiate the linear trend observed in the time-to-pattern with varying $\tau$ and the non-monotonic behavior associated with changes in domain length through numerical simulation by examining the results obtained from two distinct types of initial conditions as explained in section \ref{section4}. A  summary of the findings is presented in \cref{fig:fig5} (c) and (d) where a comparable trend is seen in the system with respect to both initial conditions. It is fascinating how the trend persists even when we opt for random initial conditions (\cref{fig:fig5} (d)). While there is inherent randomness in individual behaviors, averaging them reveals a similar linear trend with the $\tau$ and non-linear trend relating to domain size. We can draw significant conclusions regarding the variable $\mathcal{T}$. It becomes apparent that when we introduce the initial perturbation as either a random variable or to start solutions from near the unstable manifold, the time required to achieve a stable pattern increases. Under these conditions, the system exhibits a prolonged period to achieve a stable solution, thereby highlighting the system's tendency to require more time for convergence.

\begin{figure}[H]
	\begin{center}
		\includegraphics[width=0.8 \linewidth]{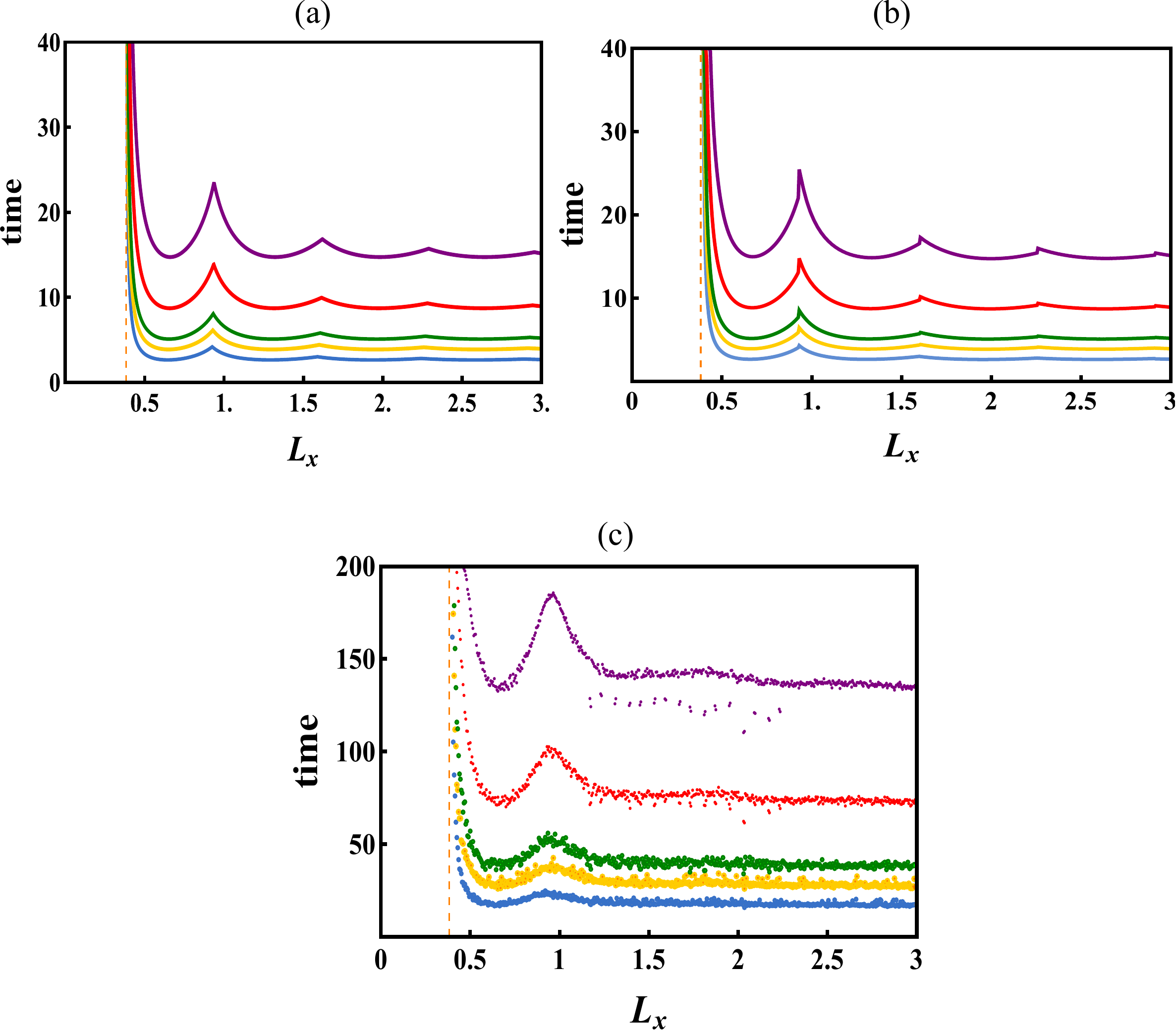}\\
		\caption{(a) Predicted time to pattern computed from the characteristic equation \eqref{eq:ch1} for $\tau = 0$ (navy blue curve), $0.1$ (yellow curve), $0.2$ (olive green curve), $0.5$ (red curve), and $1.0$ (purple curve). (b) Simulated time to pattern for the initial condition in \eqref{initial_1} when $\tau = 0$ (navy blue curve), $0.1$ (yellow curve), $0.2$ (olive green curve), $0.5$ (red curve), and $1.0$ (purple curve). (c) Simulated time to pattern from a uniform random initial condition \eqref{initial_2} for $\tau = 0$ (navy blue curve), $0.1$ (yellow curve), $0.2$ (olive green curve), $0.5$ (red curve), and $1.0$ (purple curve). Each curve in (c) is constructed using average data obtained after simulation. The orange dotted line in (a), (b), and (c) corresponds to the bifurcation point. Mode switching occurs at $L_{x} \approx 0.923$. $a=0.1$, $b =0.9$, $d_{u} = 0.01$, and $d_{v} = 0.2$. }\label{fig:fig6}
	\end{center}
\end{figure}

To conduct an in-depth analysis of the non-monotonic behavior, we have generated a time-to-pattern plot using the characteristic equation, which is visually represented in \cref{fig:fig6} (a). Each curve distinctly illustrates the non-monotonic trend with respect to domain size for various values of $\tau$. It is also evident from these curves that the duration necessary for pattern formation extends as the delay parameter is increased.

To validate these findings, we carried out supplementary numerical simulations, incorporating two distinct categories of initial conditions, as elaborated in Section \ref{section4}. \cref{fig:fig6} (b) and (c) shows the supplementary empirical support for the non-monotonic behaviour in the pattern formation dynamics in response to alterations in domain size. As shown, a similar trend is followed analogous to the trend drawn from the usage of characteristic equation. The most interesting case can be envisaged while using the random initial conditions which delineates a pattern when the random values are averaged out. 

It is also worth emphasizing that the non-monotonous behavior is more pronounced with reduced domain lengths. As the domain size expands, the prominence of this trend appears to diminish.

\subsubsection{Comparison of time to pattern for two different initial condition}
As previously mentioned, for our purpose of simulating the system we introduced two types of initial conditions. From the result shown in \cref{fig:fig5}, it is observed that the time required for the system after initiation to reach a pattern significantly differs for these two initial conditions \ref{initial_1}  and  \ref{initial_2}. Specifically, the system exhibits a notably prolonged transition time when initialized with random conditions (\ref{initial_1}) compared to when initiated from points on the unstable manifold of the equilibrium $(u_*, v_*)$. This observation implies that the system, when subjected to random initial conditions, experiences a prolonged period to depart from the equilibrium state.

When the initial conditions are situated along the unstable manifold, they are proximate to the unstable modes of the system. Consequently, trajectories corresponding to these initial conditions will swiftly diverge from the equilibrium point, as depicted in \cref{fig:fig5}(c).

Alternatively, when initial conditions are randomly selected, certain trajectories may inadvertently align with unstable modes while others do not. Trajectories closely associated with unstable modes swiftly depart from the equilibrium point, while those farther from such modes may exhibit a delayed departure, thereby extending the transition period. This phenomenon contributes to an overall lengthier transition time, as illustrated in \cref{fig:fig5} (d).

 \begin{figure}[H]
	\begin{center}
		\includegraphics[width=0.8 \linewidth]{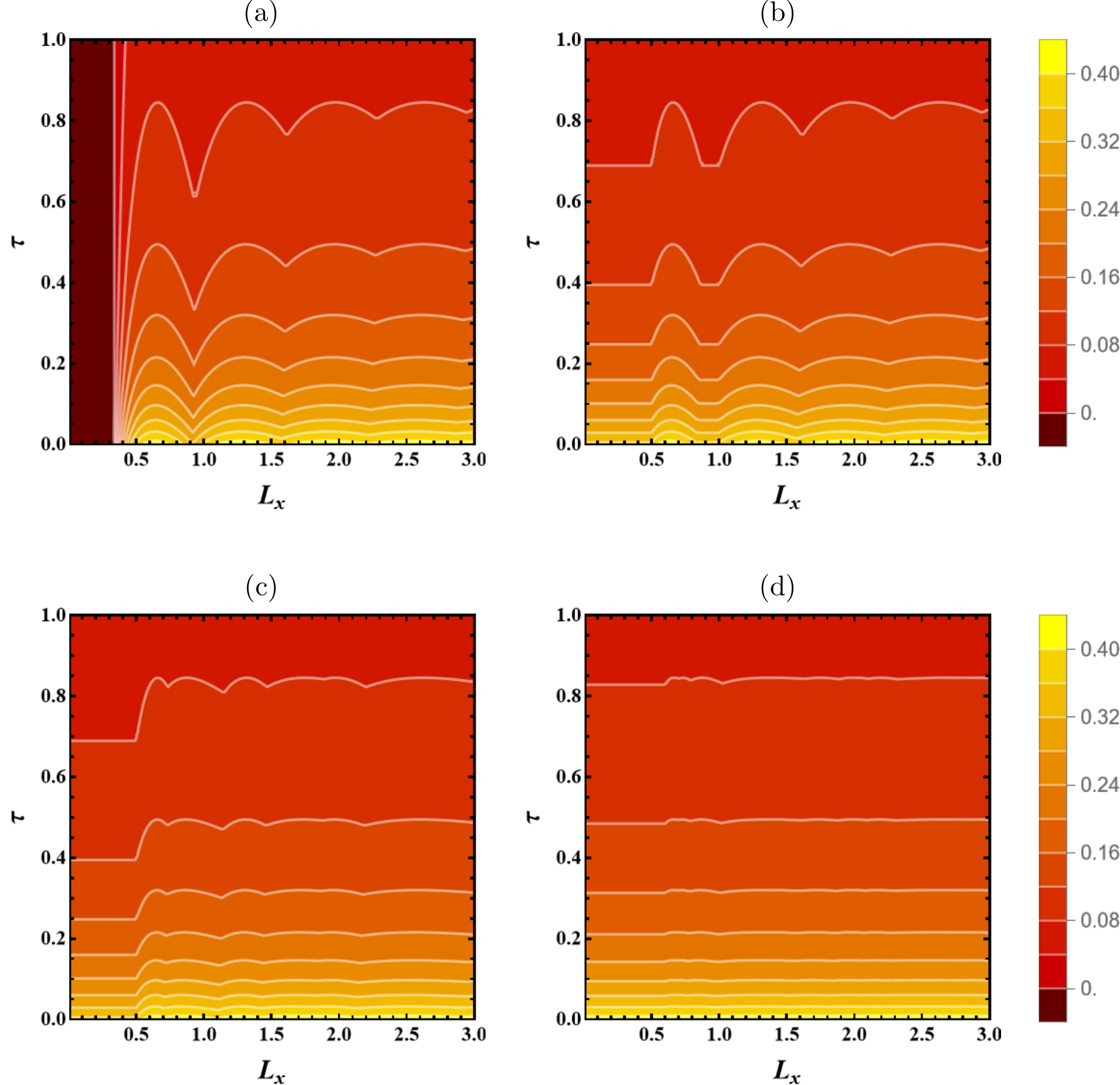}
		\caption{Heat map constructed by evaluating the $\alpha$ from the characteristic equation for the 2D LI model. $L_{x}$ is varied from $0.0$ to $3.0$ and $\tau$ from $0$ to $1$, $a = 0.1$, $b = 0.9$, $d_{u} = 0.01$, $d_{v} = 0.2$, (a) $L_{y} = 0.2$, (b)  $L_{y} = 0.5$, (c) $L_{y} = 1$ and (d) $L_{y} = 3$.}\label{fig:fig7}
	\end{center}
\end{figure}

Through \cref{fig:fig7}, we provide a detailed illustration substantiating the aforementioned statement in the form of a thermal map. This map effectively conveys the parameter $\alpha$ as defined in \eqref{eq:alpha1}, derived from the characteristic equation \eqref{eq:ch1} that depends on domain length ($L_{x}$) and temporal delay. The parameter $L_{x}$ is systematically explored across the interval of $0$ to $3$, consistently across all exhibited plots and $\tau$ from $0$ to $1$. Meanwhile, the parameter $L_{y}$ exhibits inherent variability across the four distinct plots, progressively escalating from the initial depiction in \cref{fig:fig7} (a), where $L_{y} = 0.2$, to the final instance in \cref{fig:fig7} (d), where $L_{y} = 3$.

A discernible pattern emerges from our empirical observations – the non-monotonic tendencies that are perceptible within the system's behavior (Subfigures \cref{fig:fig7} (a, b, c)) undergo a gradual process of attenuation as the dimensions of encompassing domain expands (Subfigure \cref{fig:fig7} (d)). This empirical trend underscores the significant influence of domain size on the emergent characteristics of the system, manifesting a transformative impact on the behavior elucidated within the confines of characteristic equation and its resultant $\alpha$.

\section{Onset of patterning in fixed 2D delay model}
\label{pattern_formation}

In this section, we demonstrate the occurrence of 2D pattern formation within the fixed delay LI model, while systematically varying both domain size and $\tau$. We validate the findings presented in \cref{fig:fig1} (a) and (b) by conducting numerical simulations where we identified two specific parameter combinations, denoted as $(a, b)$, which are $(0.4, 0.4)$, and $(0.1, 0.9)$. When the parameter values are set at $(a, b, \tau) = (0.4, 0.4, 0)$, linear theory predicts the absence of pattern formation. Alternatively, when the parameters are adjusted to $(a, b, \tau) = (0.1, 0.9, 0)$, the system exhibits a Turing instability, resulting in the emergence of distinct patterns. We confirm this outcome with the help of numerical simulation.

\subsection{Effects of delay}

\begin{figure}[]
	\begin{center}
		\includegraphics[width=0.8 \linewidth]{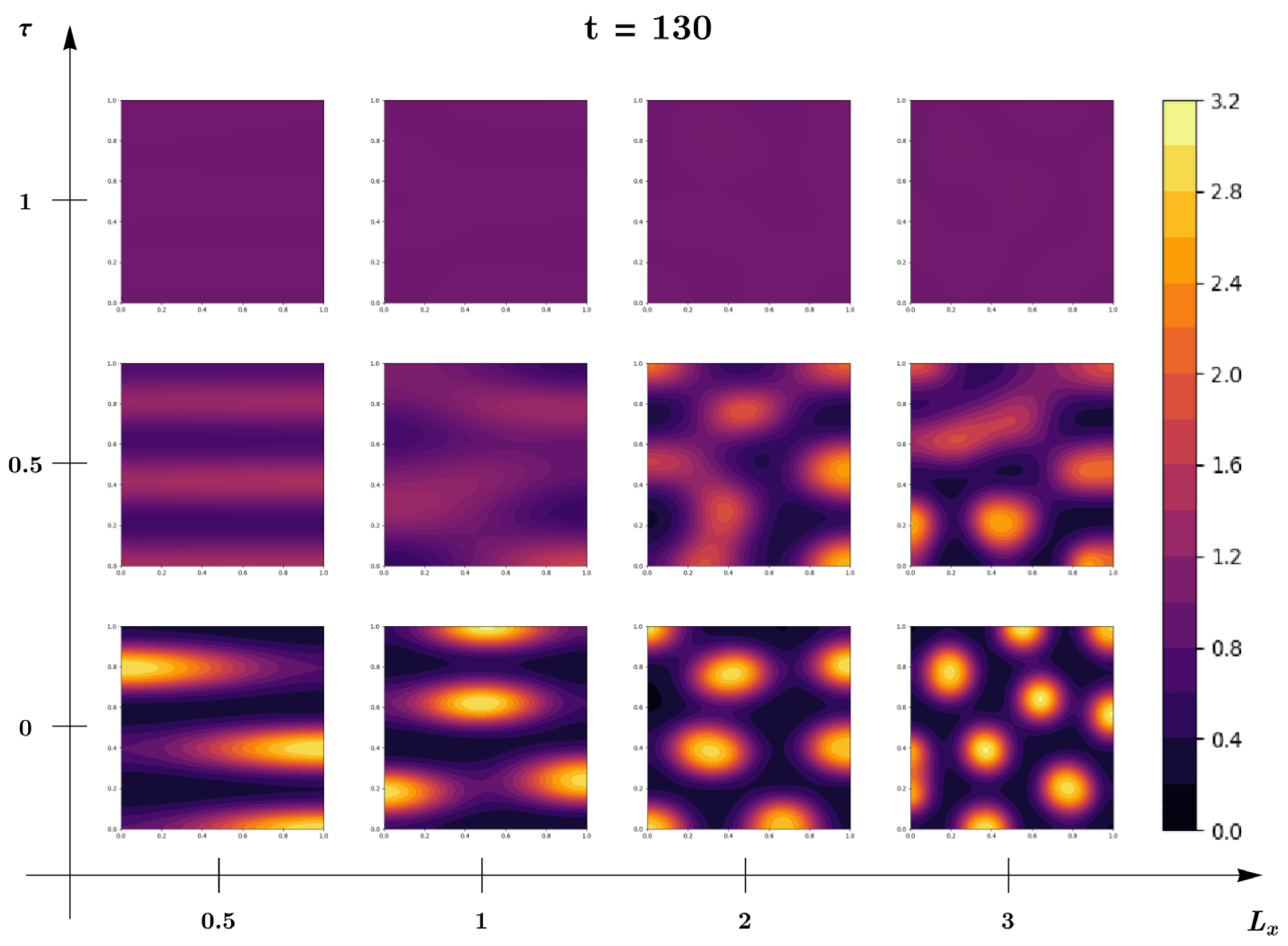}\\
		\caption{Illustrations for the pattern formation with fixed delay and domain sizes. Each snapshot in the figure corresponds to a specific temporal instance at $130$ time units. The model maintains constant parameter values: $a = 0.1$, $b = 0.9$, $d_{u} = 0.01$, and $d_{v} = 0.2$. The three sets of figure shows generating patterns for three different values of $\tau = 0, 0.5, 1$  respectively starting from the base. The effect of varying domain size is reflected laterally. $L_{y}$ if fixed to $3$, $L_{x}$ takes the values $0.5, 1, 2, 3$ respectively from the left to the right. }\label{fig:fig8a}
	\end{center}
\end{figure}

\begin{figure}[]
	\begin{center}
		\includegraphics[width=0.8 \linewidth]{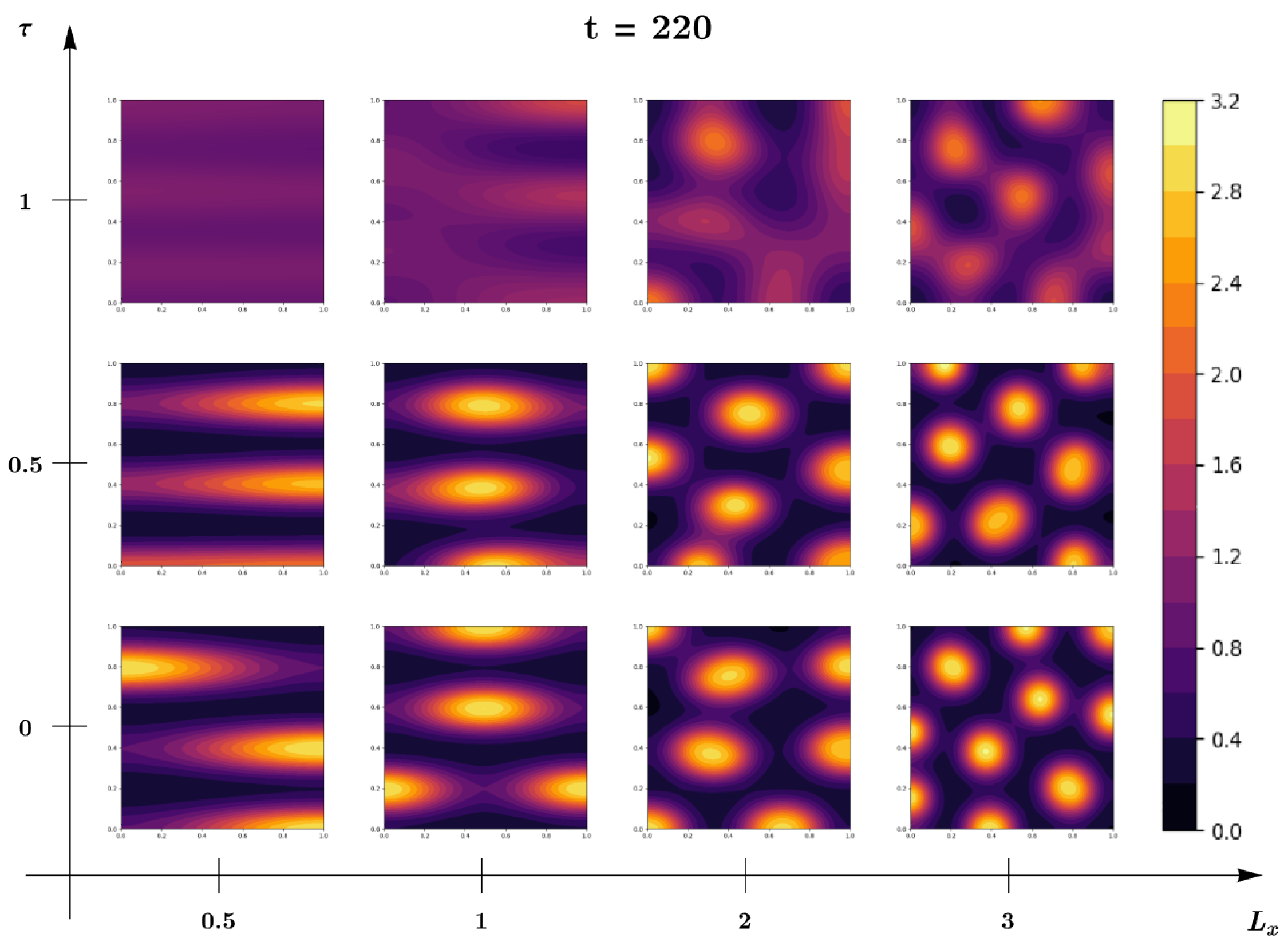}\\
		\caption{Illustrations for the pattern formation with fixed delay and domain sizes. Each snapshot in the figure corresponds to a specific temporal instance at $220$ time units. The model maintains constant parameter values: $a = 0.1$, $b = 0.9$, $d_{u} = 0.01$, and $d_{v} = 0.2$. The three sets of figure shows generating patterns for three different values of $\tau = 0, 0.5, 1$  respectively starting from the base. The effect of varying domain size is reflected laterally. $L_{y}$ if fixed to $3$, $L_{x}$ takes the values $0.5, 1, 2, 3$ respectively from the left to the right.}\label{fig:fig8b}
	\end{center}
\end{figure}

Firstly, we consider \cref{fig:fig8a} and \cref{fig:fig8b} which are the numerical plots for the LI model incorporating a constant delay. The parameter values are $(a, b) = (0.1, 0.9)$ for $\tau$ values in the set $(0, 0.5, 1)$. These solutions are complementary to the scenario with fixed delay, conforming to the equation \eqref{initial_1} with an initial small random  perturbation, over the time interval $t \in [ - \tau, 0]$. The resulting pattern is presented in a $3 \times 4$ matrix, with changing $\tau$ indicated in the vertical direction and $L_{x}$ in the horizontal direction. The corresponding time is indicated above. This figure serves as an explicit example for Subfigure \cref{fig:fig7} (d) where $L_{y}$ is fixed to 3.  \cref{fig:fig8a} illustrates the progression of the system dynamics until 130 time units. For $\tau$ is $0$, the system up to the present already shows stable pattern formation within this brief time period. On the other hand, while, for $\tau = 0.5$, at this time point, we observe the initiation of patterns, for  $\tau = 1$, there is no observable pattern formation at all. Further development of the patterns generating in the system are portrayed in  \cref{fig:fig8b}. This figure shows the status of the system after  220 time units has elapsed,  for  $\tau = 0.5$, the patterns starts to stabilize and for $\tau = 1$ the system starts to exhibit patterns. These results clearly  underline the influence of time delay in shaping system behavior.

The significance of different domain sizes on the time taken to form patterns is largely reflected here. Consider \cref{fig:fig8a}, for $\tau = 0.5$, it shows emerging patterns at different stages simultaneously for different domain length of  $L_{x}$. As we extend the domain $L_{x}$, we observe more pronounced patterns. In other words, for larger domains, the time required for patterns to emerge decreases. We can conclude that the time for pattern formation is not only dependent on the delay, but also on the domain length. This effect that domain size has on the time for formation of patterns is amplified by larger delays. However, the delay multiplies its effect only for smaller domains.  Moreover, larger domain length $L_{x}$ are adequate for circular patterns,  whereas smaller domains are insufficient, and can only accommodate stripes. Similar outcomes are also seen for $\tau = 1$. It is important to note that in \cref{fig:fig8a} and \cref{fig:fig8b}, $L_{y}$ is held constant at 3. 

When $L_{x}$ is set to $0.5$ or smaller, the available space is insufficient to support the formation of circular patterns, ergo, the patterns are more likely to be stripes. This finding provides empirical evidence supporting the proposition that when the domain length is sufficiently short, the resultant patterns are limited to stripes. This phenomenon bears a resemblance to the patterns in the tails or legs of certain animals, where space constraints similarly dictate the predominant pattern formation \cite{murray2001mathematical}.

Each simulation commences from an initial random uniform perturbation. Therefore, while the ultimate pattern outcome, whether it manifests as circles or stripes, is contingent upon the domain length, it is essential to emphasize that the specific arrangement of these circles and stripes is intricately linked to the initial conditions. As such, it is pressing to recognize that there is no means of deterministic control over the precise configuration of these patterns.

\subsection{Effects of domain size ($L_{x}$,$L_{y}$) and parameter (a, b)}

\begin{figure}[]
	\begin{center}
		\includegraphics[width=0.8 \linewidth]{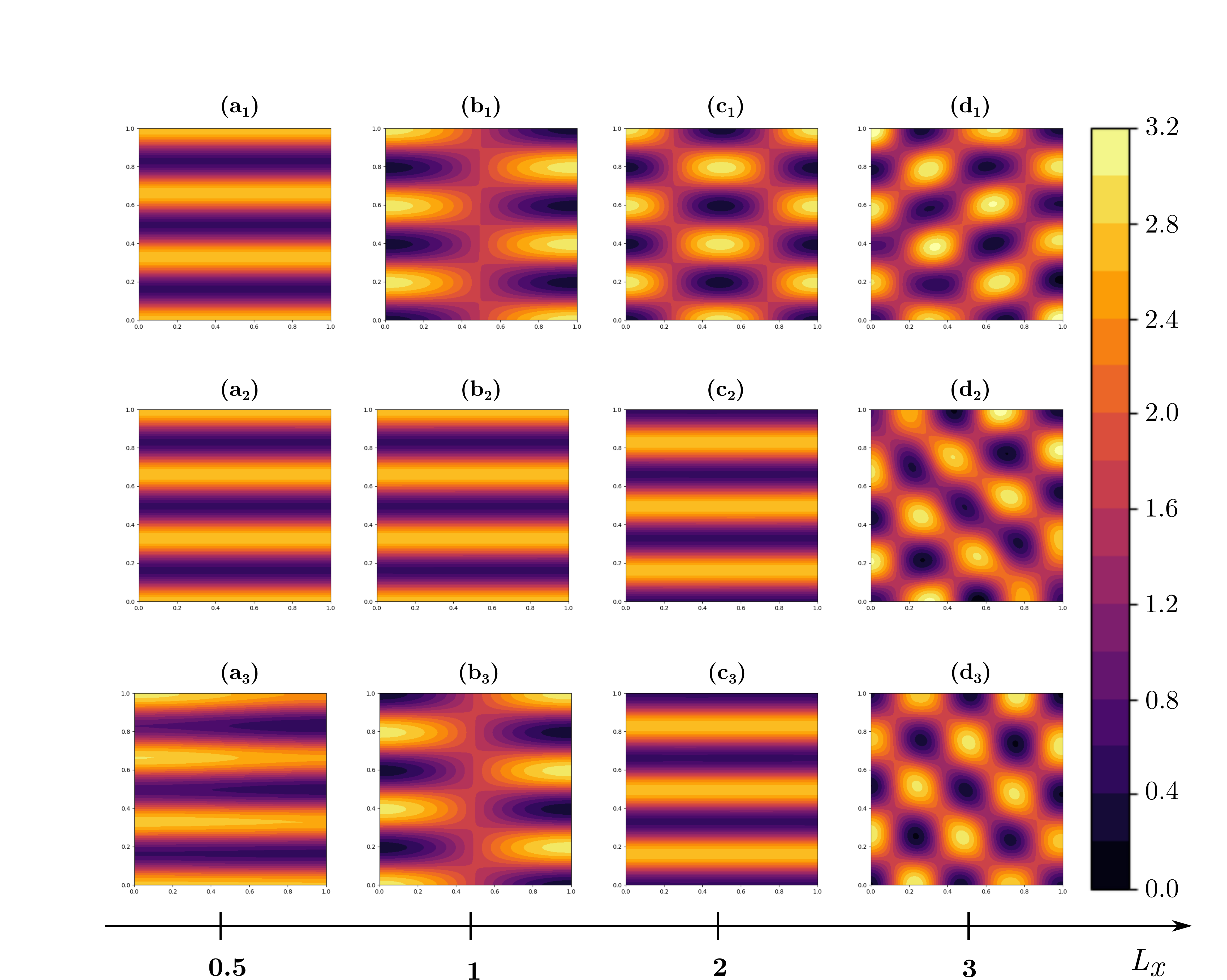}\\
		\caption{Sequential snapshots depicting pattern formation within the ligand internalisation (LI) model, employing consistent parameter settings of $a = 0.1$, $b = 1.5$, $d_{u} = 0.01$, and $d_{v} = 0.2$ showing the effects of domain and parameter $(a, b)$. T = 1000 time units for each snapshot. Subfigures ($a_1$) through ($d_3$) represent scenarios with $\tau = 0$, highlighting the absence of temporal delay. Notably, each snapshot corresponds to a distinct spatial domain size, denoted in the figure. The images in each set form a cohesive grouping, emphasizing common characteristics or evolving patterns.}\label{fig:fig_10}
	\end{center}
\end{figure}

Simulations in this section are  subject to parameter values that yield eigenvalues significantly closer to zero. In particular, we consider $(a, b) = (0.1, 1.5)$, the results of which are shown in \cref{fig:fig_10}. This figure displays distinct final solutions from three simulation runs for each of four different values of $L_X = 0.5, 1, 2, 3$, resulting in a $3 \times 4$ grid of figures. Firstly, note that, in this case, the time for a pattern to stabilize is about $1000$ time units. This observation can be attributed to $(a, b)$ values, which leads $\alpha$ values close to zero. Secondly, this set of $(a, b)$ gives rise to more complex and interesting patterns as witnessed in \cref{fig:fig_10}. The emergence of honeycomb-like structures, lines, and half circles is contingent upon factors such as the initial random mixture, domain dimensions, and the temporal evolution of the system. Interestingly, when $a$ and $b$ are set to $0.1$ and $0.9$, respectively, these factors seem to have a lesser influence. For instance, when the domain dimensions are set to $(L_{x}, L_{y}) = (2, 3)$, with $(a, b) = (0.1, 0.9)$, circular patterns appear at random positions. When $(a, b) = (0.1, 1.5)$ under the same domain conditions, stripes emerge, \cref{fig:fig_10} (c, f, g), which is unexpected given the domain size, which is large enough to equip circles, but instead exhibits stripes. This renders it a unique and noteworthy case. Even with $(L_{x}, L_{y}) = (3, 3)$, the patterns are highly organized, resembling hexagonal shapes. Thus, when $(a, b) = (0.1, 0.9)$, once the circular pattern forms, it remains stable. When $(a, b) = (0.1, 1.5)$, the system self-organizes into a more symmetric pattern over time. Thus, it can be deduced that in this particular parameter configuration, the patterns exhibit a higher degree of organization.

Another intriguing observation arising from this dataset is that the ultimate pattern outcome may exhibit variability depending on how the solution from the random initial function evolves over time. For instance, when we choose $L_x = 1$ and perform three distinct simulations, the final solutions yields different outcomes for the same $L_x$. This result becomes more interesting when $L_x = 2$. While this variability can be seen as an idiosyncrasy of the utilization of random initial solutions, it remains intriguing as one would anticipate consistent results even with the use of random perturbations. 

This observation suggests that the values of parameters $(a, b)$, along with the length of the spatial domain $(L_{x}, L_{y})$, exert a significant influence during the process of pattern formation. Notably, the final pattern outcome becomes particularly intriguing when the spatial domain length is relatively small.

\section{Discussion}
\label{discussion}
A two-component LI model is considered in a bounded 2-D domain to study the effects of discrete time delays and domain length on Turing pattern formation. One of our primary objectives is to establish how including gene expression time delays in a higher dimension can significantly impact the dynamics of a reaction-diffusion system. To this end, we carry out a methodical study of the linear stability theory and numerical simulations to investigate the behavior of pattern formation, with a focus on the influence of different time delays, shifting domain sizes, as well as variations in the initial conditions. The findings of our study enhance the understanding of the effects of delays, in relation to various other important aspects of the dynamical system in consideration. This helps build knowledge upon critical insights for both qualitative and quantitative understanding of gene regulation processes. 

As suggested in the paper by Gaffney \cite{gaffney2006gene}, when examining patterning processes that require rapid pattern establishment, such as in embryonic tissue development, the importance of considering time delays is important and highlights the need for careful justification when neglecting them. The importance of the role played by time delays in reaction kinetics has been underscored by numerous prior studies, which have consistently shown that even a minor time delay can significantly extend the length of time taken for the emergence of patterns.  The main emphasis of our study has been placed on the situation where time delay is present solely in the reaction kinetics, which is the focal point of the LI model.

We begin with the hypothesis that a linear relationship exists between the gene expression time delay and the time duration for pattern formation in a 2D LI model, a phenomenon that has been observed in the one-dimensional setting \cite{sargood2022fixed}.  Our findings indicate that the decelerating process of pattern formation due to time delay exhibits a linear scaling relationship in two-dimensions as well, which is evidenced via both linear stability analysis, and numerical investigation of the reaction-diffusion systems. 

The scope of the present study extends beyond establishing the relationship between time delays and the duration of pattern formation in a 2-D domain.  We further investigate the behavior of pattern formation processes in 2-D domains of different sizes. It is observed that there exists a non-monotonic relationship between the time required for pattern formation and the domain length, indicating the dependence of pattern formation on the defined space is a critical factor in validating the pattern-forming mechanism. The non-monotonic behavior suggests that specific spatial dimensions are more conducive to pattern formation, while others may inhibit or delay pattern emergence. These findings emphasize the significance of considering domain size as a key parameter when investigating pattern formation processes. By understanding the role that spatial extent plays, we can gain valuable insights into the underlying mechanisms that govern the development and maintenance of patterns in complex systems. In this regard, critical domain sizes that optimize the efficiency of pattern formation are identified. 

We also test a variety of starting conditions, from random disturbances of the uniform steady state to carefully chosen states from its unstable manifold. These choices can strongly influence the patterning and help us better understand the biological mechanisms behind cellular organization.

\section*{Acknowledgments}
We would like to acknowledge the assistance of Dr. Géza Makay with the code. N.~P.~D.~was supported by research grant KKP 129877, I.~Bal\'azs was supported by research grant TKP2021–NVA–09, B.~D. was supported by research grant RRF–2.3.1–21–2022–00006, G.~R. was supported by research grant TKP2021-EGA-05.

\bibliographystyle{unsrt} 
\bibliography{references}

@article{aragon1998spatial,
	title={Spatial patterning in modified {T}uring systems: Application to pigmentation patterns on marine fish},
	author={Arag{\'o}n, JL and Varea, C and Barrio, RA and Maini, PK},
	journal={Forma},
	volume={13},
	number={3},
	pages={213--221},
	year={1998}
}

@article{arcuri1986pattern,
	title={Pattern sensitivity to boundary and initial conditions in reaction-diffusion models},
	author={Arcuri, P and Murray, James D},
	doi = {10.1007/BF00275996},
	journal={Journal of Mathematical Biology},
	volume={24},
	number={2},
	pages={141--165},
	year={1986},
	publisher={Springer}
}

@article{doelman1997pattern,
	title={Pattern formation in the one-dimensional Gray-Scott model},
	author={Doelman, Arjen and Kaper, Tasso J and Zegeling, Paul A},
	doi = {10.1088/0951-7715/10/2/013},
	journal={Nonlinearity},
	number={2},
	volume={10},
	pages={523},
	year={1997},
	publisher={IOP Publishing}
}

@article{bard1981model,
	doi = {10.1016/0022-5193(81)90109-0},
	title={A model for generating aspects of zebra and other mammalian coat patterns},
	author={Bard, Jonathan BL},
	journal={Journal of Theoretical Biology},
	volume={93},
	number={2},
	pages={363--385},
	year={1981},
	publisher={Elsevier}
}

@article{barrio1999two,
	title={A two-dimensional numerical study of spatial pattern formation in interacting {T}uring systems},
	author={Barrio, RA and Varea, C and Arag{\'o}n, JL and Maini, PK},
	doi = {10.1006/bulm.1998.0093},
	journal={Bulletin of Mathematical Biology},
	volume={61},
	number={3},
	pages={483--505},
	year={1999},
	publisher={Elsevier}
}

@article{benson1993diffusion,
	title={Diffusion driven instability in an inhomogeneous domain},
	author={Benson, Debbie L and Sherratt, Jonathan A and Maini, Philip K},
	doi = {10.1016/S0092-8240(05)80270-8},
	journal={Bulletin of Mathematical Biology},
	volume={55},
	number={2},
	pages={365--384},
	year={1993},
	publisher={Elsevier}
}

@article{chen2013global,
	
	title={Global attractivity of equilibrium in Gierer--Meinhardt system with activator production saturation and gene expression time delays},
	author={Chen, Shanshan and Shi, Junping},
	doi = {10.1016/j.nonrwa.2012.12.004},
	journal={Nonlinear Analysis: Real World Applications},
	volume={14},
	number={4},
	pages={1871--1886},
	year={2013},
	publisher={Elsevier}
}

@article{chen2011note,
	title={A note on {H}opf bifurcations in a delayed diffusive {L}otka--{V}olterra predator--prey system},
	author={Chen, Shanshan and Shi, Junping and Wei, Junjie},
	doi = {10.1016/j.camwa.2011.07.011},
	journal={Computers \& Mathematics with Applications},
	volume={62},
	number={5},
	pages={2240--2245},
	year={2011},
	publisher={Elsevier}
}

@article{chen2013effect,
	title={The effect of delay on a diffusive predator-prey system with {H}olling type-II predator functional response},
	author={Chen, Shanshan and Shi, Junping and Wei, Junjie},
	doi = {10.3934/cpaa.2013.12.481},
	journal={Communications on Pure \& Applied Analysis},
	volume={12},
	number={1},
	pages={481},
	year={2013},
	publisher={American Institute of Mathematical Sciences}
}

@article{chen2013time,
	title={Time delay-induced instabilities and {H}opf bifurcations in general reaction--diffusion systems},
	author={Chen, Shanshan and Shi, Junping and Wei, Junjie},
	doi = {10.1007/s00332-012-9138-1},
	journal={Journal of Nonlinear Science},
	volume={23},
	number={1},
	pages={1--38},
	year={2013},
	publisher={Springer}
}

@article{turing1952chemical,
	title={The Chemical Basis of Morphogenesis},
	author={Alan Turing},
	doi = {10.1098/rstb.1952.0012},
	journal={Philos Trans R Soc Lond Ser B Biol Sci},
	number={},
	volume={237},
	year={1952},
	month={},
	pages={37–72},
}

@article{gierer1972theory,
	doi = {10.1007/BF00289234},
	title={A theory of biological pattern formation},
	author={Gierer, Alfred and Meinhardt, Hans},
	journal={Kybernetik},
	volume={12},
	number={1},
	pages={30--39},
	year={1972},
	publisher={Springer}
}

@article{murray1981pattern,
	title={On pattern formation mechanisms for lepidopteran wing patterns and mammalian coat markings},
	author={Murray, James Dickson},
	doi = {10.1098/rstb.1981.0155},
	journal={Philosophical Transactions of the Royal Society of London. B, Biological Sciences},
	volume={295},
	number={1078},
	pages={473--496},
	year={1981},
	publisher={The Royal Society London}
}

@article{murray1981pre,
	title={A pre-pattern formation mechanism for animal coat markings},
	author={Murray, JD610574},
	doi = {10.1016/0022-5193(81)90334-9},
	journal={Journal of Theoretical Biology},
	volume={88},
	number={1},
	pages={161--199},
	year={1981},
	publisher={Elsevier}
}

@article{cocho1987discrete,
	title={Discrete systems, cell-cell interactions and color pattern of animals. I. Conflicting dynamics and pattern formation.},
	author={Cocho, Germinal and P{\'e}rez-Pascual, Rafael and Rius, Jos{\'e} L},
	doi = {10.1016/S0022-5193(87)80211-4},
	journal={Journal of Theoretical Biology},
	volume={125},
	number={4},
	pages={419--435},
	year={1987}
}

@article{meinhardt1987model,
	title={A model for pattern formation on the shells of molluscs},
	author={Meinhardt, Hans and Klingler, Martin},
	doi = {10.1016/S0022-5193(87)80101-7},
	journal={Journal of Theoretical Biology},
	volume={126},
	number={1},
	pages={63--89},
	year={1987},
	publisher={Elsevier}
}

@article{nijhout1978wing,
	title={Wing pattern formation in {L}epidoptera: a model},
	author={Nijhout, H Frederik},
	doi = {10.1002/jez.1402060202},
	journal={Journal of Experimental Zoology},
	volume={206},
	number={2},
	pages={119--136},
	year={1978},
	publisher={Wiley Online Library}
}

@article{meinhardt1998models,
	title={Models of pattern formation applied to plant development},
	author={Meinhardt, Hans and Koch, Andr{\'e}-Joseph and Bernasconi, Giuliano},
	doi = {10.1142/9789814261074_0027},
	journal={Symmetry in Plants},
	volume={4},
	pages={723--758},
	year={1998},
	publisher={Singapore: World Scientific}
}

@article{izhikevich2006fitzhugh,
	title={{F}itzhugh-{N}agumo model},
	author={Izhikevich, Eugene M and FitzHugh, Richard},
	journal={Scholarpedia},
	volume={1},
	number={9},
	pages={1349},
	year={2006}
}

@article{prigogine1968symmetry,
	title={Symmetry breaking instabilities in dissipative systems. II},
	author={Prigogine, Ilya and Lefever, Ren{\'e}},
	doi = {10.1063/1.1668896},
	journal={The Journal of Chemical Physics},
	volume={48},
	number={4},
	pages={1695--1700},
	year={1968},
	publisher={American Institute of Physics}
}

@article{kondo1995reaction,
	title={A reaction--diffusion wave on the skin of the marine angelfish {P}omacanthus},
	author={Kondo, Shigeru and Asai, Rihito},
	doi = {10.1038/376765a0},
	journal={Nature},
	volume={376},
	pages={765--768},
	year={1995},
	publisher={Springer}
}

@article{plaza2004effect,
	title={The effect of growth and curvature on pattern formation},
	author={Plaza, Ram{\'o}n G and Sanchez-Garduno, Faustino and Padilla, Pablo and Barrio, Rafael A and Maini, Philip K},
	doi = {10.1007/s10884-004-7834-8},
	journal={Journal of Dynamics and Differential Equations},
	volume={16},
	pages={1093--1121},
	year={2004},
	publisher={Springer}
}

@article{kim2020pattern,
	title={Pattern formation in reaction--diffusion systems on evolving surfaces},
	author={Kim, Hyundong and Yun, Ana and Yoon, Sungha and Lee, Chaeyoung and Park, Jintae and Kim, Junseok},
	doi = {10.1016/j.camwa.2020.08.026},
	journal={Computers \& Mathematics with Applications},
	volume={80},
	number={9},
	pages={2019--2028},
	year={2020},
	publisher={Elsevier}
}

@article{jensen1994localized,
	doi = {10.1103/PhysRevE.50.736},
	title={Localized structures and front propagation in the {L}engyel-{E}pstein model},
	author={Jensen, O and Pannbacker, Viggo Ole and Mosekilde, Erik and Dewel, Guy and Borckmans, Pierre},
	journal={Physical Review E},
	volume={50},
	number={2},
	pages={736},
	year={1994},
	publisher={APS}
}

@article{schnakenberg1979simple,
	doi ={10.1016/0022-5193(79)90042-0},
	title={Simple chemical reaction systems with limit cycle behaviour},
	author={Schnakenberg, J},
	journal={Journal of Theoretical Biology},
	volume={81},
	number={3},
	pages={389--400},
	year={1979},
	publisher={Elsevier}
}

@article{toral2003characterization,
	doi = {10.1016/S0378-4371(03)00198-5},
	title={Characterization of the anticipated synchronization regime in the coupled {F}itzhugh--{N}agumo model for neurons},
	author={Toral, Ra{\'u}l and Masoller, Cristina and Mirasso, Claudio R and Ciszak, Marzena and Calvo, O},
	journal={Physica A: Statistical Mechanics and its Applications},
	volume={325},
	number={1-2},
	pages={192--198},
	year={2003},
	publisher={Elsevier}
}

@article{sargood2022fixed,
	doi = {10.1007/s11538-022-01052-0},
	title={Fixed and Distributed Gene Expression Time Delays in Reaction-Diffusion Systems},
	author={Sargood, Alec and Gaffney, Eamonn A and Krause, Andrew L},
	journal={Bulletin of Mathematical Biology},
	year={2022, 84.9: 98}
}

@article{gaffney2006gene,
	doi = {10.1007/s11538-006-9066-z},
	title={Gene expression time delays and {T}uring pattern formation systems},
	author={Gaffney, EA and Monk, NAM},
	journal={Bulletin of Mathematical Biology},
	volume={68},
	number={1},
	pages={99--130},
	year={2006},
	publisher={Springer}
}

@article{lee2010aberrant,
	title={Aberrant behaviours of reaction diffusion self-organisation models on growing domains in the presence of gene expression time delays},
	author={Lee, S Seirin and Gaffney, EA},
	doi = {10.1007/s11538-010-9533-4},
	journal={Bulletin of Mathematical Biology},
	volume={77},
	pages={2161-2179},
	year={2010}
}

@article{tennyson1995human,
	doi = {10.1038/ng0295-184},
	title={The human dystrophin gene requires 16 hours to be transcribed and is cotranscriptionally spliced},
	author={Tennyson, Christine N and Klamut, Henry J and Worton, Ronald G},
	journal={Nature Genetics},
	volume={9},
	number={2},
	pages={184--190},
	year={1995},
	publisher={Nature Publishing Group}
}

@article{rybakin2001morphogenesis,
	title={Morphogenesis and pattering in hydra. II. Molecular mechanisms},
	author={Rybakin, VS},
	journal={Tsitologiia},
	volume={43},
	number={1},
	pages={39--45},
	year={2001}
}

@article{green2002morphogen,
	doi = {10.1002/dvdy.10170},
	title={Morphogen gradients, positional information, and {X}enopus: interplay of theory and experiment},
	author={Green, Jeremy},
	journal={Developmental dynamics: an official publication of the American Association of Anatomists},
	volume={225},
	number={4},
	pages={392--408},
	year={2002},
	publisher={Wiley Online Library}
}

@article{grieneisen2012morphogengineering,
	doi = {10.1186/1752-0509-6-37},
	title={Morphogengineering roots: comparing mechanisms of morphogen gradient formation},
	author={Grieneisen, Ver{\^o}nica A and Scheres, Ben and Hogeweg, Paulien and M Mar{\'e}e, Athanasius F},
	journal={BMC Systems Biology},
	volume={6},
	number={1},
	pages={1--20},
	year={2012},
	publisher={Springer}
}

@article{hadeler2007interaction,
	doi = { 10.3934/dcdsb.2007.8.95},
	title={Interaction of diffusion and delay},
	author={Hadeler, Karl Peter and Ruan, Shigui},
	journal={Discrete \& Continuous Dynamical Systems-B},
	volume={8},
	number={1},
	pages={95},
	year={2007},
	publisher={American Institute of Mathematical Sciences}
}

@article{morita1984destabilization,
	doi = {10.1007/BF03167861},
	title={Destabilization of periodic solutions arising in delay-diffusion systems in several space dimensions},
	author={Morita, Yoshihisa},
	journal={Japan Journal of Applied Mathematics},
	volume={1},
	number={1},
	pages={39--65},
	year={1984},
	publisher={Springer}
}

@article{ruan1998turing,
	doi = {10.1093/imamat/61.1.15},
	title={{T}uring instability and travelling waves in diffusive plankton models with delayed nutrient recycling},
	author={Ruan, SHIGUI},
	journal={IMA Journal of Applied Mathematics},
	volume={61},
	number={1},
	pages={15--32},
	year={1998},
	publisher={Oxford University Press}
}

@article{ruan1999persistence,
	doi = {10.1006/jdeq.1998.3599},
	title={Persistence and extinction in two species reaction-diffusion systems with delays},
	author={Ruan, Shigui and Zhao, Xiao-Qiang},
	journal={Journal of Differential Equations},
	volume={156},
	number={1},
	pages={71--92},
	year={1999},
	publisher={Elsevier}
}

@article{gourley1996predator,
	doi = {10.1007/BF00160498},
	title={A predator-prey reaction-diffusion system with nonlocal effects},
	author={Gourley, SA and Britton, NF},
	journal={Journal of Mathematical Biology},
	volume={34},
	number={3},
	pages={297--333},
	year={1996},
	publisher={Springer}
}

@article{freedman1997global,
	doi = {10.1006/jdeq.1997.3264},
	title={Global asymptotics in some quasimonotone reaction-diffusion systems with delays},
	author={Freedman, HI and Zhao, Xiao-Qiang},
	journal={Journal of Differential Equations},
	volume={137},
	number={2},
	pages={340--362},
	year={1997},
	publisher={Elsevier}
}

@article{higham1995existence,
	doi = {10.1016/0168-9274(95)00051-U},
	title={Existence and stability of fixed points for a discretised nonlinear reaction-diffusion equation with delay},
	author={Higham, Desmond J and Sardar, Tasneem},
	journal={Applied Numerical Mathematics},
	volume={18},
	number={1-3},
	pages={155--173},
	year={1995},
	publisher={Elsevier}
}

@article{yi2017bifurcation,
	doi = {10.3934/dcdsb.2017031},
	title={The bifurcation analysis of {T}uring pattern formation induced by delay and diffusion in the {S}chnakenberg system},
	author={Yi, Fengqi and Gaffney, Eamonn A and Seirin-Lee, Sungrim},
	journal={Discrete \& Continuous Dynamical Systems-B},
	volume={22},
	number={2},
	pages={647},
	year={2017},
	publisher={American Institute of Mathematical Sciences}
}

@book{murray2001mathematical,
	title={Mathematical Biology II: Spatial Models and Biomedical Applications computational model predicts phenotype from genotype},
  author={Murray, James D},
  journal={(No Title)},
  year={2003}
}

@article{lewis2003autoinhibition,
	doi = {10.1016/S0960-9822(03)00534-7},
	title={Autoinhibition with transcriptional delay: a simple mechanism for the zebrafish somitogenesis oscillator},
	author={Lewis, Julian},
	journal={Current Biology},
	volume={13},
	number={16},
	pages={1398--1408},
	year={2003},
	publisher={Elsevier}
}

@article{jing2006mechanisms,
	doi = {10.1016/j.mod.2006.03.006},
	title={Mechanisms underlying long-and short-range nodal signaling in Zebrafish},
	author={Jing, Xiao-hong and Zhou, Sheng-mei and Wang, Wei-qing and Chen, Yu},
	journal={Mechanisms of Development},
	volume={123},
	number={5},
	pages={388--394},
	year={2006},
	publisher={Elsevier}
}

@article{sander2003pattern,
doi = {10.1016/S0022-0396(02)00156-0},
  title={Pattern formation in a nonlinear model for animal coats},
  author={Sander, Evelyn and Wanner, Thomas},
  journal={Journal of Differential Equations},
  volume={191},
  number={1},
  pages={143--174},
  year={2003},
  publisher={Elsevier}
}
\end{document}